\documentclass[11pt]{article}

\usepackage{bm}
\usepackage{booktabs}
\usepackage{algorithm}  
\usepackage{changepage}
\usepackage{xcolor}
\usepackage{algorithmic}

\usepackage{epsfig}

\usepackage{amsmath}
\usepackage{amsfonts}
\usepackage{amssymb}
\usepackage{dsfont}
\usepackage{mathrsfs}
\usepackage{bbm}
\usepackage{hyperref}
\usepackage{amsmath}
\usepackage{amsfonts}
\usepackage{amssymb}
\usepackage{multirow}
\usepackage{mathrsfs}
\usepackage{natbib}

\numberwithin{equation}{section}

\newtheorem{theorem}{Theorem}[section]
\newtheorem{lemma}{Lemma}[section]
\newtheorem{proposition}{Proposition}[section]
\newtheorem{example}{Example}[section]
\newtheorem{remark}{Remark}[section]
\newtheorem{definition}{Definition}[section]
\newtheorem{assumption}{Assumption}[section]

\newcommand{\setd}{{ d \kern -.15em l}}
\newcommand{\hatsetd}{ d \hat{\kern -.15em l }}
\newcommand{\dd}{\mathsf {d\kern -0.07em l}} 

\newcommand{\N}{{\cal N}}

\newcommand{\bgeqn}{\begin{eqnarray}}
\newcommand{\edeqn}{\end{eqnarray}}
\newcommand{\bgeq}{\begin{eqnarray*}}
\newcommand{\edeq}{\end{eqnarray*}}
\newcommand{\bec}{\begin{center}}
\newcommand{\enc}{\end{center}}
\newcommand{\R}{{\rm I\!R}}

\newcommand{\inmat}[1]{\mbox{\rm {#1}}}

\newcommand{\U}{{\cal U}}

\newcommand{\Z}{{\cal Z}}
\newcommand{\F}{{\cal F}}

\newcommand{\C}{{\cal C}}
\newcommand{\A}{{\cal A}}

\newcommand{\M}{{\cal M}}
\newcommand{\X}{{\cal X}}

\newcommand{\Y}{{\cal Y}}

\newcommand{\be}{\begin{equation}}
\newcommand{\ee}{\end{equation}}

\def\w{\omega}

\title{
Existence and Uniqueness Theorem of
Continuous and Monotone Bayesian Nash Equilibrium 
and Stability Analysis\footnote{This project is supported by a CUHK start-up grant. }
}

\author{ 
Ziheng Su\footnote{Department of Systems Engineering and Engineering Management, The Chinese University of Hong Kong, Shatin, N.T., Hong Kong.  
Email: zsu@se.cuhk.edu.hk.} \, and \,
Huifu Xu\footnote{Department of Systems Engineering and Engineering Management, The Chinese University of Hong Kong, Shatin, N.T., Hong Kong. 
Email: hfxu@se.cuhk.edu.hk.}
}

\begin{document}

\allowdisplaybreaks

\maketitle

\begin{abstract}
Since the seminal work by Meirowitz~\cite{meirowitz2003existence}, there has been growing attention on the existence and uniqueness of continuous Bayesian Nash equilibria. In the existing literature, existence is typically established using Schauder’s fixed-point theorem, relying on the equicontinuity of players’ best response functions. Uniqueness, on the other hand, is usually derived under additional monotonicity conditions. In this paper, we revisit the issues of existence and uniqueness, and advance the literature by establishing both simultaneously using the Banach fixed-point theorem under a set of moderate conditions. Furthermore, we analyze the stability of such equilibria with respect to perturbations in the joint probability distribution of type parameters, offering theoretical support for the application of Bayesian Nash equilibrium models in data-driven contexts.
\end{abstract}


\section{Introduction}

In recent years, there has been a growing focus on Nash games involving private information. A common assumption in such games is that the initial distribution of all players' types is publicly known. Each player has full knowledge of their own type, which determines their utility function, while remaining unaware of their opponents' types. Based on the prior information, each player selects a response strategy within their type space based on the Nash conjecture, and an equilibrium
arising from this kind of game is known as a Bayesian Nash equilibrium (BNE). The concept of the Bayesian game was introduced by Harsanyi~\cite{harsanyi1967_1, harsanyi1968_2, harsanyi1968_3} and has attracted considerable attention, spanning from modeling to equilibrium analysis and applications in strategic decision-making under incomplete information.

Existence of pure strategy Bayesian Nash equilibria (BNEs) has been studied extensively ever since.
A majority of research on the existence of monotone BNE (MBNE) has focused on conditions related to supermodularity and single-crossing condition (SCC).
Vives~\cite{vives1990nash} analyzes existence of Nash equilibria of non-cooperative games by using lattice theoretical methods, where each player's payoffs satisfy certain supermodularity properties which are directly related to strategic complementarities. The author also establishes existence of BNE as the supermodularity is preserved by integration.
Athey~\cite{athey2001single} studied Bayesian games with one-dimensional actions and type spaces, assuming that types are drawn from an atomless joint probability distribution. 
When each player's expected utility function exhibits SCC property, i.e., each player’s optimal response is increasing with respect to the increase of its type parameter in response to rivals’ increasing strategies (for type varying), she demonstrated existence of pure strategy Nash equilibria (PSNEs) in every finite action game with each player’s optimal response strategy being increasing and step-like. Moreover, when action spaces are continuous, Athey shows existence of a sequence of increasing step-like PSNEs of finite action games that converges to an equilibrium with the continuum-action, which means that an equilibrium in continuous action spaces can be approximated by a sequence of increasing step-like PSNEs in finite action spaces. 
McAdams~\cite{mcadams2003isotone} extends Athey’s results into more general settings by considering the quasi-supermodularity (in own action) and single-crossing (in own action and type), where both action and type spaces are considered to be multidimensional and only partially ordered.
Reny~\cite{reny2011existence} further relaxes the requirements on payoff functions, showing that monotone best-reply conditions can hold under weaker assumptions.
Van Zandt and Vives~\cite{vives2007monotone}
introduce a concept called ``monotone supermodularity'' in Bayesian games, where, apart from strategic complementarities,
each player’s payoff also has increasing differences in their own action and the profile of types, and that each player’s posterior beliefs are increasing in type with respect to first-order stochastic dominance.
They provide an induction method to construct a greatest and a least BNE, with both equilibria being monotone in type. 
Mason and Valentinyi~\cite{mason2007existence} further study monotone pure-strategy equilibria in Bayesian games and derive existence and uniqueness under a contraction-mapping framework. Their analysis is focused on the monotone equilibrium formulation.

Another research stream focuses on 
continuous BNEs where 
each player's optimal response function is continuous.
By using Schauder’s fixed-point theorem under some moderate conditions, 
Meirowitz~\cite{meirowitz2003existence} establishes existence of a continuous BNE (CBNE) with each player's equilibrium strategy being continuous. One of the main conditions is the equicontinuity of the set of optimal response strategies, which ensures that the operator mapping the set of optimal response strategies to itself is compact. Meirowitz comments that the equicontinuity condition is either likely to fail or difficult to verify in practical applications.
Ui~\cite{ui2016bayesian} provides a sufficient condition for existence and uniqueness of a differentiable BNE by regarding it as a solution of a variational inequality problem. The author demonstrates that when the Jacobian matrix of the payoff gradient is negative definite for each state, a BNE is unique.
Guo et al.~\cite{guo2021existence} take a step forward by presenting some verifiable conditions for the required equicontinuity, namely, some growth conditions of the expected utility function of each player. Moreover, under some monotonicity conditions, the authors demonstrate the uniqueness of a CBNE. The authors also develop some computational approaches for finding an approximate CBNE by using the polynomial decision rule to restrict the strategy function.

In a more recent development, the
CBNE model 
is extended to the case where each player's 
action space depends on 
the type parameter and/or rivals' actions.
For instance, Tao and Xu~\cite{tao2024generalized} consider a generalized BNE model
where each player's action space 
is represented by a system of 
inequalities involving both the player's own type parameter and the rivals' actions.
Under 
some moderate
conditions, they demonstrate
existence of continuous GBNE and the uniqueness of such an equilibrium when each player’s action space is only dependent on its type.
In the case that each player’s action space is also dependent on rivals’ actions, they give a simple example to show that the uniqueness of GBNE is not guaranteed under standard monotone conditions. 
In a further study, Su and Xu~\cite{SU2025BNRE} introduce a Bayesian Nash risk equilibrium (BNRE) model, where they 
differentiate the underlying uncertainties by 
considering aleatoric uncertainty and epistemic uncertainty, and discussing the model from two perspectives: monotonicity of a continuous BNRE and continuity of a monotone BNRE. 
This line of research has recently been extended further by Liu~\cite{liu2026games}, who proposes a risk-revising Bayesian Nash equilibrium and studies its existence and comparative statics with respect to risk aversion.


In this paper, we follow the strand of research but with different focuses. 
First, we derive 
existence and uniqueness of BNE by Banach's fixed-point theorem. 
Unlike Schauder's fixed-point theorem and Tarski's fixed-point theorem, Banach's fixed-point theorem relies on the contraction property of the optimal response operator, which allows us to establish existence and uniqueness simultaneously.  
Mason and Valentinyi~\cite{mason2007existence} and Mathevet~\cite{mathevet2010contraction} seem to be the first
to use the contraction-based argument to derive existence and uniqueness of a BNE albeit the argument has 
been widely used in other game theoretic and equilibrium models, see, e.g., \cite{van2000existence,gao2017price,graf2024stochastic}.
Second, we investigate the stability of BNEs when the probability distribution of type parameters is perturbed for various reasons.
As far as we are concerned, the paper makes the following contributions.


\begin{itemize}

\item \textbf{Existence and uniqueness 
of BNE.}
We use the Banach fixed-point theorem to derive three versions of existence and uniqueness results
under the strong concavity of each player's  utility function (Assumption~\ref{ass:equi}):
BNE with measurable optimal response functions (Theorem~\ref{thm:BNE-contraction}), BNE with Lipschitz continuous optimal response functions (Theorem~\ref{thm:exist-unique-cbne}) and BNE with monotone optimal response functions (Theorem~\ref{thm:exist-unique-CMBNE}).
The sufficient condition for uniqueness complements 
the well-known strict monotonicity condition Ui~\cite[Lemma 4 and Proposition 2]{ui2016bayesian} (Example~\ref{ex:asym_2p_BG}). 
Moreover, unlike the existing research in the literature~\cite{guo2021existence,tao2024generalized,SU2025BNRE} which utilizes growth condition of player's objective function to derive H\"older continuity of the optimal response function, 
we use the implicit function theorem of a variational system 
to derive Lipschitz continuity of each player's optimal response function with respect to both their own type parameter and the actions of their rivals (Theorem~\ref{thm:lipschitz-theta}).
The Lipschitz continuity provides a potential new avenue 
for computing an approximate CBNE and paves the way for 
the follow-up stability analysis.

 \item  \textbf{Stability of BNE.}
In practical application of BNE models, the true probability distribution of type parameters is often unknown particularly in data-driven problems. This raises a question as to whether
a BNE derived from a theoretical model based on the true unknown distribution is sufficiently close to the one based on 
a perturbed/approximate distribution. 
In this paper, we 
derive an error bound of the BNE based on the true and 
perturbed distributions via stability analysis
of stochastic variational inequality systems
(Theorems~\ref{thm:belief-perturb} and \ref{thm:belief-perturb-2}). 
These new results 
fill an important gap in the literature of 
BNE and provide theoretical 
guarantees for applying BNE models to data-driven problems, where the probability distribution 
in theoretical models (often treated as training models) is shifted.


\end{itemize}

The rest of the paper is organized as follows. 
In Section~\ref{se:The model}, we introduce the Bayesian game and define the continuous and monotone Bayesian Nash equilibrium (CMBNE) concept. 
Section~\ref{se:Lipschitz} establishes the Lipschitz continuity of the optimal response function.
Section~\ref{se:Existence and uniqueness} presents our main existence and uniqueness results for BNE and discusses how our conditions relate to the existing literature. 
Section~\ref{se:Monotone} studies the monotonicity of the equilibrium.
In Section~\ref{se:stability}, we study the stability of CBNE with respect to perturbations of the joint distribution of types. 
Finally, Section~\ref{se:Concluding remarks} concludes the paper and outlines directions for future research.

Throughout the paper, we adopt the following notation. By convention,
we use $\R^n$ to denote $n$-dimensional Euclidean space, and  
$\R_+^n$ and $\R_{++}^n$ to represent
the first  orthant and its interior  respectively.
In the case that $n=1$, we drop the superscript.
 We write
$\|x\|$ and $\|x\|_\infty$ to denote the Euclidean norm and the infinity norm of vector $x \in \R^n$ respectively and 
the same notation for the induced matrix norms.
Let $(X, \mathscr{B},\mu)$ be a measurable space and $p\in [1,+\infty]$. We use
$
{\cal L}^{p}(X,\mu)
:=
\left\{
f : X \to \R^n \ \text{measurable}
\;:\;
\int_X \|f(x)\|^{p} \, d\mu(x) < \infty
\right\} 
$
to denote the ${\cal L}^{p}$-space of functions mapping from $(X, \mathscr{B})$ to $\R^n$.
For each $f\in {\cal L}^{p}(X,\mu)$, we 
write 
$
\|f\|_{{\cal L}^{p}(\mu)}
:=
\left(
\int_X \|f(x)\|^{p} \, d\mu(x)
\right)^{1/p},
$ 
and 
$
\|f\|_{{\cal L}^\infty(\mu)}
:=
\operatorname*{ess\,sup}_{x \in X} \|f(x)\|
$
for the ${\cal L}^{p}(\mu)$-norm and ${\cal L}^\infty(\mu)$-norm. 
Let
$w\in\R_{++}^n$ and $\mu:=\Pi_{i=1}^n \mu_i$.
For $f=(f_1,\cdots,f_n)$, where $f_i$ is a vector-valued measurable function in ${\cal L}^{p}(X_i,\mu_i)$, 
we use  
$
\|f\|_{\infty,w,{\cal L}^{p}(\mu)}
:=
\max_{1\leq i\leq n}
\frac{1}{w_i}\left(
\|f_i\|_{{\cal L}^{p}(\mu_i)}
\right)$
to denote the weighted ${\cal L}^{p}(\mu)$-norm
of $f$.
We also use $\| f_i\|_{\infty} = \sup_{\theta_i \in \Theta_i} \|f_i(\theta_i)\|$
to denote the infinity norm
and $\| f \|_{\infty,w}
:=
\max_{1\leq i\leq n}
\left(
\frac{1}{w_i}\| f_i \|_\infty
\right)$
to denote the weighted infinity norm.
In the case that $f_i$ is continuous, $\| f_i \|_{\infty}=
\|f_i\|_{{\cal L}^\infty(\mu_i)}$.
Finally, we use  $\mathscr{P}(X)$ to denote the set of all probability measures over 
$(X,\mathscr{B})$.


\section{The Model}
\label{se:The model}

We consider a Bayesian game with $n$ players.
Each player possesses a preference utility function
denoted by $u_i(a_i, a_{-i}, \theta)$, for $i \in N :=\{1,\cdots,n\}$,
which depends on the player's
action $a_i$, its rival's actions $a_{-i}:= (a_1, \cdots, a_{i-1}, a_{i+1}, \cdots, a_n)$, and
the player $i$'s type $\theta_i$.
We assume that a type $\theta_i$ takes values from set $\Theta_i$
and an action $a_i$ takes values from action space $\A_i$,
where $\Theta_i$ and $\A_i$ are non-empty, compact and convex subsets
of $\R^{d_i}$ and $\R^{z_i}$ respectively.
Following the terminology of Meirowitz~\cite{meirowitz2003existence}, a profile of types is a vector $\theta := (\theta_1,\cdots,\theta_n)\in \Theta := \Theta_1\times \cdots\times \Theta_n$ and a profile of actions is a vector $a := (a_1,\cdots,a_n)\in \A := \A_1\times \cdots\times \A_n$.
Using the standard notation, we denote by $a_{-i}$ and $\theta_{-i}$ respectively the vector of actions and the types of all players except $i$.

Following the standard setup in Bayesian games, we assume that $\theta$ is a random vector endowed with a prior distribution $\eta \in \mathscr{P}(\Theta)$
over 
full support (i.e., $\operatorname{supp}(\eta) = \Theta$). The prior $\eta$ is common knowledge among all players.
This information describes the probability of all players taking a particular $\theta$, which may be retrieved from empirical data in practice.
For each $i \in N$, let $\eta_i \in \mathscr{P}(\Theta_i)$ denote the marginal distribution of $\theta_i$.
Information on players' types is private, which means each player only knows its own type but not others. 
After observing its own type $\theta_i$, player $i$'s posterior belief about $\theta_{-i}$ is represented by a conditional probability distribution $\eta_i(\cdot \mid \theta_i) \in \mathscr{P}(\Theta_{-i})$, which describes the probability of player $i$'s rivals taking a particular type $\theta_{-i}$.
Throughout the paper, we will use  $\theta$ to denote 
both a deterministic element of $\R^{d_1+\cdots+d_n}$ 
and a random vector $\theta(\w): (\Omega,\mathscr{B},\eta) \rightarrow \R^{d_1+\cdots+d_n}$ depending on the context.

In Bayesian games with pure strategy, the response (strategy) function of player $i$ is represented by $f_i$ mapping from type space $\Theta_i$ to action space $\A_i$.
For $i \in N$, we denote by $\F_i$ the set of measurable functions $f_i : \Theta_i \to \A_i$ endowed with the infinity norm $\| f_i\|_{\infty} = \sup_{\theta_i \in \Theta_i} \|f_i(\theta_i)\|$, and $\C_i \subset \F_i$  
the set of continuous functions,
where $\|\cdot\|$ denotes the Euclidean norm.
Moreover, we use 
$\F_i^+ \subset \F_i$ and 
$\C_i^+ \subset \C_i$
to denote the sets of increasing functions and increasing continuous functions respectively,
and we will call functions in $\F^+/\C^+$ also simply monotone.
For the simplicity of notation,
we write
\begin{subequations}
\begin{align}
& \F \
(\inmat{resp.}~\F^+) :=\prod_{i\in N} \F_i\
(\inmat{resp.}~\F_i^+), \quad
\F _{-i}\
(\inmat{resp.}~\F_{-i}^+)
:=
\prod_{j \in N\backslash \{i\}} 
\F_j\
(\inmat{resp.}~\F_j^+),
\label{eq:F}
\\
& \C \
(\inmat{resp.}~\C^+) 
:=\prod_{i\in N}
\C_i\
(\inmat{resp.}~\C_i^+),
\quad 
\C _{-i}\
(\inmat{resp.}~\C_{-i}^+)
:= 
\prod_{j \in N\backslash \{i\}}
\C_j\
(\inmat{resp.}~\C_j^+).
\label{eq:F-C}
\end{align}
\end{subequations} 

We define each player’s expected utility function given the rivals' strategy function $f_{-i}$ by
\begin{equation}
\label{eq:def-phi_i-BNE}
\phi_i(a_i, f_{-i}, \theta_i) := \int_{\Theta_{-i}} u_i(a_i, f_{-i}(\theta_{-i}), \theta_i, \theta_{-i}) d \eta_i \left(\theta_{-i} | \theta_i\right), \
\forall \theta_i \in \Theta_i,
\end{equation}
where the expectation is taken with respect to the conditional probability distribution $\eta_i \left(\theta_{-i} | \theta_i\right)$, 
and $f_{-i}(\theta_{-i})$ represents the rivals’ action at scenario 
$\theta_{-i}\in \Theta_{-i}:=
\Theta_1\times \cdots \times \Theta_{i-1}\times\Theta_{i+1}\times \cdots \times \Theta_n$.  
Assuming that each player chooses its optimal strategy by maximizing expected utility under the Nash conjecture (taking rivals' optimal response functions  as given), we consider a situation in which no player can benefit from unilaterally deviating from its chosen strategy. This leads to the formal definition of a Bayesian Nash equilibrium.

\begin{definition}[Bayesian Nash equilibrium (BNE)]
\label{def:BNE}
\rm A
{\em  pure strategy Bayesian Nash equilibrium}
is an $n$-tuple of functions $f^*:= (f^*_1, \cdots, f^*_n)\in \F$
mapping from $\Theta_1\times\cdots\times \Theta_n$ to
$\A_1\times \cdots \times \A_n$
such that
\begin{equation}
\label{eq:BNE}
\inmat{(BNE)}  \quad f^*_i(\theta_i) \in\mathop{\arg \max}_{a_i\in \A_i}
\int_{\Theta_{-i}} u_i(a_i, f^*_{-i}(\theta_{-i}),\theta_i, \theta_{-i}) d\eta_i
(\theta_{-i}|\theta_i), \ \forall \theta_i\in \Theta_i,
\end{equation}
for $i\in N$. 
$f^*$ is said to be a 
continuous and monotone BNE (CMBNE)
if $f_i^*$ is continuous and monotone over $\Theta_i$ for $i\in N$. 
\end{definition}

Based on \eqref{eq:def-phi_i-BNE}, we can write the BNE model \eqref{eq:BNE} succinctly as 
\begin{align}
\label{eq:BNE2}
f^*_i(\theta_i) \in \mathop{\arg \max}_{a_i\in \A_i}\  
\phi_i(a_i, f^*_{-i},\theta_i), 
\ \forall \theta_i\in \Theta_i,
\end{align}
for $ i \in N$.
Model (\ref{eq:BNE2}) can be 
interpreted purely from Nash equilibrium perspective, 
where each player's optimal action 
is assumed to be a function of 
player's type parameter, rather than a fixed action in a Nash game.
This is an infinite-dimensional 
Nash equilibrium problem where each player's optimal action is a function. 
In the forthcoming discussions, we will derive sufficient conditions for existence and uniqueness of three types of 
BNEs where (a) $f^*$ is merely measurable (Section~\ref{se:Existence and uniqueness}), (b) $f^*$ is Lipschitz continuous (Section~\ref{se:Existence and uniqueness}), and (c) $f^*$ is monotonically increasing (Section~\ref{se:Monotone}).

\section{Lipschitz continuity of optimal response function}
\label{se:Lipschitz}

Unlike the existing literature, we establish the existence and uniqueness of BNEs primarily by exploiting the Lipschitz continuity of each player’s optimal response function with respect to variations in the rivals’ optimal response functions. Moreover, we also require the Lipschitz continuity with respect to each player’s own type parameter to derive existence of a Lipschitz continuous BNE. In this section, we provide the technical details on the two types of Lipschitz continuity.
To this end, we make the following assumption.

\begin{assumption}
\label{assu:u_i}
Consider the BNE model \eqref{eq:BNE}.
For $i \in N$, 
(a) $u_i(a, \theta)$ is continuous in $(a, \theta)$, 
$\nabla_{a_i} u_i(a, \theta)$ is continuous in $a$ and measurable in $\theta$; 
(b) for each $\hat{a}_i\in\A_i$,
there exist a neighborhood $\U_i(\hat{a}_i;f_{-i},\theta_i)$ of $\hat{a}_i$  relative to 
$\A_i$, which 
depends on $f_{-i}$ and $\theta_i$,
and a nonnegative function
$h_{\hat{a}_i,f_{-i}}:\Theta_{i} \times \Theta_{-i} \to \R_+$
with
$h_{\hat{a}_i,f_{-i}}(\theta_i, \cdot \ ) \in \mathcal{L}^1
\big( \eta_i(\cdot \mid \theta_i) \big)
$
such that
\begin{equation}
\sup_{a_i\in \U_i(\hat{a}_i;f_{-i},\theta_i)}
\big\|\nabla_{a_i}u_i(a_i,f_{-i}(\theta_{-i}),\theta_i,\theta_{-i})\big\|
\leq
h_{\hat{a}_i,f_{-i}}(\theta_i, \theta_{-i}),
\ \forall \theta_{-i} \in \Theta_{-i};
\end{equation}
(c) for every Borel set $S \in \mathscr{B}(\Theta_{-i})$, the conditional distribution $\eta_i(S \mid \cdot)$ is measurable on $\Theta_i$.
\end{assumption}

Let $\phi_i(a_i,f_{-i},\theta_i)$
be defined as in \eqref{eq:def-phi_i-BNE}. 
Under Assumption~\ref{assu:u_i}, 
$\phi_i(a_i,f_{-i},\theta_i)$ is continuous in $a_i$ over $\A_i$
for each fixed $f_{-i}\in\F_{-i}$ and $\theta_i\in\Theta_i$, and it is  
measurable in $\theta_i$ over $\Theta_i$ for each fixed $f_{-i}\in\F_{-i}$ and $a_i\in\A_i$  (see Bertsekas~\cite[Proposition A.4]{bertsekas2012dynamic}). 
This means  
that $\phi_i$ is a Carath\'eodory function of $(a_i,\theta_i)$.
By the measurable maximum theorem~\cite[Theorem~18.19]{aliprantis2006infinite},
for each $f_{-i}\in \F_{-i}$, the correspondence
$\theta_i \mapsto \mathop{\arg\max}_{a_i\in\A_i} \phi_i(a_i,f_{-i},\theta_i)$ is measurable with nonempty compact values, and therefore admits a measurable selector.
Moreover, for each fixed $f_{-i}\in\F_{-i}$ and $\theta_i\in\Theta_i$,
 we can 
 use the dominated convergence theorem to show, 
  under
 Assumption~\ref{assu:u_i} (b), that
 $\phi_i(a_i, f_{-i}, \theta_i)$
is differentiable on $\A_i$ with
\begin{equation}
\label{eq:gradient_phi_ai}
\nabla_{a_i}
\phi_i(a_i, f_{-i}, \theta_i)
=
\int_{\Theta_{-i}} \nabla_{a_i}u_i(a_i,f_{-i}(\theta_{-i}),\theta_i,\theta_{-i}) d\eta_i(\theta_{-i}|\theta_i).
\end{equation} 
A stronger but often easier-to-verify condition that implies the same conclusion is the global Lipschitz continuity in $a_i$, i.e., there exists $h_i:\Theta_i \times \Theta_{-i} \to \R_+$ with $h_i(\theta_i, \cdot\ ) \in \mathcal{L}^1 \big(\eta_i(\cdot\mid\theta_i)\big)$ such that, for all $a_i',a_i''\in\A_i$ and all $\theta_{-i}\in\Theta_{-i}$,
\begin{equation}
\big| u_i(a_i',f_{-i}(\theta_{-i}),\theta_i,\theta_{-i}) 
- u_i(a_i'',f_{-i} (\theta_{-i}),\theta_i,\theta_{-i}) \big| 
\leq
h_i(\theta_i,\theta_{-i}) \|a_i'-a_i''\|.
\end{equation}
Together with Assumption~\ref{assu:u_i}(a), this implies 
$\|\nabla_{a_i}u_i(a_i,f_{-i}(\theta_{-i}),\theta_i,\theta_{-i})\|
\leq h_i(\theta_i, \theta_{-i})$ for all $a_i \in \A_i$, so Assumption~\ref{assu:u_i}(b) holds with $h_{\hat{a}_i, f_{-i}}\equiv h_i$.
In the forthcoming discussions, we use Assumption~\ref{assu:u_i}.

\subsection{Lipschitz continuity in rivals' strategies}
To obtain a global Lipschitz estimate for the player's optimal response function, we further make the following assumption.

\begin{assumption}
\label{ass:equi}
Consider the BNE model (\ref{eq:BNE}). 
For $i \in N$, assume:
(a) $u_i(a_i, a_{-i}, \theta)$ is strongly concave over $\A_i$ uniformly for all $a_{-i}$, i.e., there exists a positive constant
$\sigma_i$ such that
\begin{equation}
\label{eq:strongly-concave}
\left(
\nabla_{a_i} u_i(a_i', a_{-i}, \theta)  
-
\nabla_{a_i} u_i(a_i'', a_{-i}, \theta) 
\right)^\top 
\left(a_i' - a_i''\right)
\leq
-\sigma_i \|a_i' - a_i'' \|^2,\
\forall a_i', a_i'' \in \A_i;
\end{equation}
(b) $\nabla_{a_i} u_i(a_i, a_{-i}, \theta)$ is 
Lipschitz continuous  in $a_{-i}$ over $\A_{-i}$  uniformly for all $a_i$ and $\theta$ with moduli $\{ \tau_{ij} \}_{j \neq i}$, i.e.,
\begin{equation}
\label{eq:lipschitz-2}
\|
\nabla_{a_i}u_i(a_i, a_{-i}', \theta) - 
\nabla_{a_i}u_i(a_i, a_{-i}'', \theta)
\|   \leq
\sum_{j \neq i} \tau_{ij}\|a_{j}' - a_{j}'' \|,\ \forall a_{-i}', a_{-i}'' \in \A_{-i}.
\end{equation}
\end{assumption}
The constant $\sigma_i$ in inequality \eqref{eq:strongly-concave} represents 
the minimal curvature of the utility function, which signifies the ratio of 
the marginal utility w.r.t.~player $i$'s change of action.  
The condition implies that the utility function $u_i$ is strongly concave in $a_i$ uniformly w.r.t.~$a_{-i}$ and $\theta$.
The strong concavity condition
and the compactness of ${\cal A}_i$ ensure 
existence and uniqueness of the optimal solution 
\begin{equation}
\label{eq:optimal-response}
\A_i^*(f_{-i}, \theta_i)
:= \mathop{\arg \max}_{a_i\in \A_i}
\int_{\Theta_{-i}} u_i(a_i, f_{-i}(\theta_{-i}),\theta_i, \theta_{-i}) \, d\eta_i(\theta_{-i}|\theta_i).
\end{equation}
Moreover, it follows from Berge's maximum theorem (see e.g.,~\cite[Theorem 17.31]{aliprantis2006infinite}) that
$\A_i^*(f_{-i}, \theta_i)$ is continuous in $f_{-i}$ under the infinity norm.  Under Assumption~\ref{assu:u_i},
$\A_i^*(f_{-i}, \theta_i)$ is measurable in $\theta_i$, that is, $\A_i^*(f_{-i},\cdot\ ) \in \F_i$.
For each fixed $f_{-i}\in\F_{-i}$, we call the mapping $\theta_i \mapsto \A_i^*(f_{-i},\theta_i)$ player $i$'s {\em optimal response function}.
Assumption~\ref{ass:equi}(b) 
imposes a blockwise Lipschitz bound on the marginal utility $\nabla_{a_i}u_i$ with moduli
$\{\tau_{ij}\}_{j\ne i}$, controlling the cross–effect of rivals’ actions on player $i$’s marginal utility.
Together with (a), this linear control yields a global (blockwise) Lipschitz estimate for player $i$’s
optimal response with respect to rivals’ strategies $f_{-i}$. 
The next theorem addresses this.

\begin{theorem}[Lipschitz continuity w.r.t. rivals' strategies]
\label{thm:lipschitz-f}
Consider the BNE model \eqref{eq:BNE}.
Suppose that Assumptions~\ref{assu:u_i} - \ref{ass:equi} hold.
Then the optimal response function $\A_i^*(f_{-i},\cdot)\in\F_i$ is blockwise Lipschitz continuous in $f_{-i}$ with moduli $\{\frac{\tau_{ij}}{\sigma_i}\}_{j\neq i}$, 
in the following sense:
\begin{itemize}
    \item[(i)]for any $p\in[1, +\infty]$,
\begin{equation}
\label{eq:lipschitz-f}
\|\A_i^*(f_{-i},\cdot)-\A_i^*(g_{-i},\cdot)\|_{\mathcal L^p(\eta_i)}
\leq
\sum_{j\neq i}\frac{\tau_{ij}}{\sigma_i}
\|f_j-g_j\|_{\mathcal L^p(\eta_j)},
\quad \forall f_{-i},g_{-i}\in\F_{-i};
\end{equation}

    \item[(ii)]under the infinity norm,
\begin{equation}
\label{eq:lipschitz-f2}
\|\A_i^*(f_{-i},\cdot)-\A_i^*(g_{-i},\cdot)\|_\infty
\leq
\sum_{j\neq i}\frac{\tau_{ij}}{\sigma_i}
\|f_j-g_j\|_\infty,
\quad \forall f_{-i},g_{-i}\in\F_{-i}.
\end{equation}
\end{itemize} 
\end{theorem}

\noindent
\textbf{Proof.}
We use Lemma~\ref{Lem:Lip-slv-mon} (in the Appendix) to prove the result.
Under Assumptions~\ref{assu:u_i} - \ref{ass:equi}, we can
integrate both sides of \eqref{eq:strongly-concave} with respect to $\eta_i(\cdot\mid\theta_i)$ and use \eqref{eq:gradient_phi_ai} to
 obtain that $-\nabla_{a_i}\phi_i(\cdot,f_{-i},\theta_i)$
is strongly monotone over $\A_i$ with modulus $\sigma_i$, uniformly for all $f_{-i}$ and $\theta_i$. 
Moreover, for each fixed $f_{-i}\in\F_{-i}$ and $\theta_i\in\Theta_i$, $\A_i^*(f_{-i},\theta_i)$ satisfies the following first-order optimality condition:
\begin{equation}
\label{eq:phi-vi}
0 \in
-\nabla_{a_i} \phi_i \big( \A_i^*(f_{-i}, \theta_i), f_{-i}, \theta_i \big)
+ \N_{\A_i} \big(\A_i^*(f_{-i},\theta_i) \big),
\end{equation}
where $\N_{\A_i}$ denotes the convex normal cone to ${\cal A}_i$.
Therefore, by fixing $\theta_i$ and treating 
$f_{-i}$ as a varying parameter, we can apply
Lemma~\ref{Lem:Lip-slv-mon} to establish 
stability of $-\nabla_{a_i} \phi_i\left(\cdot, f_{-i}, \theta_i\right)$ w.r.t. the variation of $f_{-i}$ as follows: 
\begin{align}
\nonumber
& \|\A_i^*(f_{-i}, \theta_i) - \A_i^*(g_{-i}, \theta_i) \| \\
\nonumber
& \leq 
\frac{1}{\sigma_i }
\left \|
\nabla_{a_i} \phi_i\left(\A_i^*(f_{-i}, \theta_i), f_{-i}, \theta_i\right)
-
\nabla_{a_i} \phi_i(\A_i^*(f_{-i}, \theta_i), g_{-i}, \theta_i)
\right \| \\ 
\nonumber
& =
\frac{1}{\sigma_i }
\left \|
\int_{\Theta_{-i}}
\big(
\nabla_{a_i} u_i\left(\A_i^*(f_{-i}, \theta_i), f_{-i}(\theta_{-i}), \theta \right)
-
\nabla_{a_i} u_i\left(\A_i^*(f_{-i}, \theta_i), g_{-i}(\theta_{-i}), \theta \right)
\big)
d \eta_i(\theta_{-i} | \theta_i)
\right \| \\ 
& \leq
\frac{1}{\sigma_i}
\int_{\Theta_{-i}}
\sum_{j \neq i}
\tau_{ij}
\|f_{j}(\theta_{j}) - g_{j}(\theta_{j})\|
d \eta_i(\theta_{-i} | \theta_i),
\label{eq:Lip-resp-rival}
\end{align}
where the last inequality is due to the Lipschitz continuity in \eqref{eq:lipschitz-2}.
Consequently, for $p\in [1,\infty)$, we have,
\begin{align}
\nonumber
& \| \A_i^*(f_{-i}, \cdot) - \A_i^*(g_{-i}, \cdot)  \|_{\mathcal{L}^p(\eta_i)} 
= 
\left(
\int_{\Theta_i} 
\|\A_i^*(f_{-i}, \theta_i) - \A_i^*(g_{-i}, \theta_i) \|^p 
d\eta_i(\theta_i) \right)^{1/p} \\ \nonumber
& \leq
\left(
\int_{\Theta_i}
\left(
\int_{\Theta_{-i}}
\sum_{j \neq i}
\frac{\tau_{ij}}{\sigma_i}
\|f_{j}(\theta_{j}) - g_{j}(\theta_{j})\|
d \eta_i(\theta_{-i} | \theta_i) \right)^p 
d \eta_i(\theta_i) \right)^{1/p} 
\quad (\inmat{by \eqref{eq:Lip-resp-rival}})
\\ \nonumber
& =
\left(
\int_{\Theta_i}
\left(
\sum_{j \neq i}
\int_{\Theta_{-i}}
\frac{\tau_{ij}}{\sigma_i}
\|f_{j}(\theta_{j}) - g_{j}(\theta_{j})\|
d \eta_i(\theta_{-i} | \theta_i) \right)^p 
d \eta_i(\theta_i) \right)^{1/p} \\ \nonumber
& \leq
\sum_{j \neq i}
\frac{\tau_{ij}}{\sigma_i}
\left(
\int_{\Theta_i}
\left(
\int_{\Theta_{-i}}
\|f_{j}(\theta_{j}) - g_{j}(\theta_{j})\|
d \eta_i(\theta_{-i} | \theta_i) \right)^p 
d \eta_i(\theta_i) \right)^{1/p} 
\quad (\inmat{by Minkowski inequality})
\\ \nonumber
& \leq
\sum_{j \neq i}
\frac{\tau_{ij}}{\sigma_i}
\left(
\int_{\Theta_i}
\int_{\Theta_{-i}}
\|f_{j}(\theta_{j}) - g_{j}(\theta_{j})\| ^p
d \eta_i(\theta_{-i} | \theta_i)  
d \eta_i(\theta_i) \right)^{1/p} 
\quad (\inmat{by Jensen's inequality})
\nonumber\\
& = 
\sum_{j \neq i}
\frac{\tau_{ij}}{\sigma_i}
\left(
\int_\Theta
\|f_{j}(\theta_{j}) - g_{j}(\theta_{j})\| ^p
d \eta(\theta) \right)^{1/p}
=
\sum_{j \neq i}
\frac{\tau_{ij}}{\sigma_i}
\| f_j - g_j\|_{\mathcal{L}^p(\eta_j)},
\end{align}
which is \eqref{eq:lipschitz-f}.
For the case $p=+\infty$, we can take the essential supremum on both sides of  \eqref{eq:Lip-resp-rival} to obtain
$$
\|\A_i^*(f_{-i},\cdot)-\A_i^*(g_{-i},\cdot)\|_{\mathcal L^\infty(\eta_i)}
\leq
\sum_{j\neq i}\frac{\tau_{ij}}{\sigma_i}
\|f_j-g_j\|_{\mathcal L^\infty(\eta_j)}.
$$
Likewise, we 
can take
the pointwise supremum in \eqref{eq:Lip-resp-rival} to establish \eqref{eq:lipschitz-f2}.
$\hfill \square$

\noindent
Theorem~\ref{thm:lipschitz-f} gives rise to two notions of Lipschitz continuity of the optimal response function under two different norms. 
It is well known that the ${\cal L}^p(\eta_i)$-norm, for $p\in[1,\infty)$, measures the “average” (or area-based) distance between the graphs of two functions, whereas for $p=\infty$, the ${\cal L}^\infty(\eta_i)$-norm is the essential supremum norm and still identifies functions that differ only on sets of $\eta_i$-measure zero.
In addition, \eqref{eq:lipschitz-f2} gives a stronger estimate under the infinity norm, which controls the maximal deviation uniformly over all $\theta_i\in\Theta_i$.
Lipschitz continuity with respect to the former allows inequality~\eqref{eq:Lip-resp-rival} to fail on a set of $\theta_i$ with zero $\eta_i$-measure. 
By contrast, Lipschitz continuity under the infinity norm requires the inequality to hold uniformly for every $\theta_i\in\Theta_i$. This distinction leads to two different notions of BNE, which we will discuss later in Section~\ref{se:Existence and uniqueness}. We will also use the Lipschitz
continuity to discuss stability of BNEs in Section~\ref{se:stability}.

\subsection{Lipschitz continuity in own type}

We now turn to discuss Lipschitz continuity of 
optimal response function with respect to each player's own type parameters. To this end, 
we introduce a new set of conditions  on 
$u_i$ and $\eta_i(\cdot \mid \theta_i)$  in the next assumption.

\begin{assumption}
\label{ass:lipschitz}
Consider the BNE model~\eqref{eq:BNE}. For each $i\in N$, assume:
(a) $\sup_{a\in \A, \theta \in \Theta}
\left \| \nabla_{a_i} u_i(a, \theta) \right\| < +\infty$;
(b) there exists a nonnegative measurable function 
$\nu_i:\Theta_{-i}\to \R_+$ such that for all $a\in\A$, 
$\theta_i',\theta_i''\in\Theta_i$, and $\theta_{-i}\in\Theta_{-i}$,
$$
\|\nabla_{a_i}u_i(a,\theta_i',\theta_{-i})
-
\nabla_{a_i}u_i(a,\theta_i'',\theta_{-i})\|
\leq
\nu_i(\theta_{-i})\|\theta_i'-\theta_i''\|,
$$
and
$$
\overline{\nu}_i:=
\sup_{\theta_i\in\Theta_i}
\int_{\Theta_{-i}} \nu_i(\theta_{-i})
\,d\eta_i(\theta_{-i} | \theta_i) < +\infty;
$$
(c) the conditional distribution $\eta_i(\cdot\mid\theta_i)$ is Lipschitz continuous in $\theta_i$ under the total variation (TV) distance with modulus $\gamma_i>0$, that is,
$$
\dd_{\mathrm{TV}}
\big(\eta_i(\cdot\mid\theta_i'),
\eta_i(\cdot\mid\theta_i'')\big)
\leq
\gamma_i\|\theta_i'-\theta_i''\|,
\quad \forall \theta_i',\theta_i''\in\Theta_i,
$$
where the TV distance between two probability measures $\mu,\nu \in \mathscr P(\Xi)$ is
$$
\dd_{\mathrm{TV}}(\mu,\nu) :=
\frac{1}{2} \sup\left\{
\left|\int_\Xi h\,d\mu-\int_\Xi h\,d\nu\right|
:\ h:\Xi\to \R\ \text{measurable},\ \|h\|_\infty\leq 1
\right\}.
$$ 
\end{assumption}

\begin{proposition}
\label{prop:lipschitz}
Consider the BNE model (\ref{eq:BNE}). 
Let Assumptions \ref{assu:u_i} and \ref{ass:lipschitz} hold.
Then $\nabla_{a_i} \phi_i(a_i, f_{-i}, \theta_i)$ is Lipschitz continuous over $\Theta_i$ for all $a_i$ and $f_{-i}$ with modulus $\kappa_i$, i.e.,
\begin{equation}
\label{eq:lipschitz-1}
\|
\nabla_{a_i} \phi_i(a_i, f_{-i}, \theta_i') - 
\nabla_{a_i} \phi_i(a_i, f_{-i}, \theta_i'')
\|   \leq
\kappa_i \|\theta_i' - \theta_i'' \|,\ \forall \theta_i', \theta_i'' \in \Theta_i,
\end{equation}    
where 
$\kappa_i := 
\overline{\nu}_i
+
2\gamma_i \sup\limits_{a\in \A, \theta \in \Theta}
\left \| \nabla_{a_i} u_i(a, \theta) \right\| $.
\end{proposition}

\noindent
\textbf{Proof.}
By Assumption \ref{assu:u_i}, we have 
equation~\eqref{eq:gradient_phi_ai}, that is,
$$
\nabla_{a_i} \phi_i(a_i, f_{-i}, \theta_i)
= \int_{\Theta_{-i}}
\nabla_{a_i} u_i(a_i, f_{-i}(\theta_{-i}), \theta_i, \theta_{-i})
d \eta_i(\theta_{-i} | \theta_i).
$$
For any $\theta_i', \theta_i'' \in \Theta_i$,
\begin{align*}
& \|
\nabla_{a_i} \phi_i(a_i, f_{-i}, \theta_i') - 
\nabla_{a_i} \phi_i(a_i, f_{-i}, \theta_i'')
\| \\
& =
\left\| \int_{\Theta_{-i}}
\nabla_{a_i} u_i(a_i, f_{-i}(\theta_{-i}), \theta_i', \theta_{-i})
d \eta_i(\theta_{-i} | \theta_i')
-
\int_{\Theta_{-i}}
\nabla_{a_i} u_i(a_i, f_{-i}(\theta_{-i}), \theta_i'', \theta_{-i})
d \eta_i(\theta_{-i} | \theta_i'') \right\| \\
& \leq
\int_{\Theta_{-i}}
\left \|
\nabla_{a_i} u_i(a_i, f_{-i}(\theta_{-i}), \theta_i', \theta_{-i})
-
\nabla_{a_i} u_i(a_i, f_{-i}(\theta_{-i}), \theta_i'', \theta_{-i})
\right\|
d \eta_i(\theta_{-i} | \theta_i') \\
& \quad +
\left\| \int_{\Theta_{-i}}
\nabla_{a_i} u_i(a_i, f_{-i}(\theta_{-i}), \theta_i'', \theta_{-i})
d \eta_i(\theta_{-i} | \theta_i')
-
\int_{\Theta_{-i}}
\nabla_{a_i} u_i(a_i, f_{-i}(\theta_{-i}), \theta_i'', \theta_{-i})
d \eta_i(\theta_{-i} | \theta_i'') \right\|
 \\
&  \leq
\left( \int_{\Theta_{-i}}
\nu_i(\theta_{-i})
d \eta_i(\theta_{-i} | \theta_i')
\right)
\left \|
\theta_i' - \theta_i''
\right\|
+
2 \dd_{\mathrm{TV}}
\big(\eta_i(\cdot \mid\theta_i'),
\eta_i(\cdot \mid\theta_i'')\big)
\left( \sup_{a_\in \A, \theta \in \Theta}
\left \|  \nabla_{a_i} u_i(a, \theta) \right\| \right) 
 \\
&  \leq
\left(
\overline{\nu}_i 
+
2\gamma_i
\sup_{a_\in \A, \theta \in \Theta}
\left \|  \nabla_{a_i} u_i(a, \theta) \right\|
\right)
\left \|
\theta_i' - \theta_i''
\right\|,
\end{align*}
which implies the Lipschitz continuity of 
$\nabla_{a_i} \phi_i(a_i, f_{-i}, \theta_i)$ in $\theta_i$
with modulus $\kappa_i$.
$ \hfill \square$

Proposition \ref{prop:lipschitz} shows that the marginal expected utility
$\nabla_{a_i}\phi_i(a_i,f_{-i},\theta_i)$ is Lipschitz in the type $\theta_i$ (modulus $\kappa_i$),
providing stability in the own–type dimension and thereby bounding how the player’s marginal
expected utility varies with $\theta_i$. 
The next theorem shows the Lipschitz continuity of the optimal response function $\A_i^*(f_{-i}, \theta_i)$ with respect to $\theta_i$ .

\begin{theorem}[Lipschitz continuity w.r.t. own type]
\label{thm:lipschitz-theta}
Consider the BNE model \eqref{eq:BNE}.
Suppose that Assumptions~\ref{assu:u_i}, \ref{ass:equi}(a) and \ref{ass:lipschitz} hold,
then for each $f_{-i} \in \F_{-i}$, the optimal response $\A_i^*(f_{-i}, \theta_i)$ is Lipschitz continuous in $\theta_i$ with modulus $\frac{\kappa_i}{\sigma_i}$, that is,
\begin{equation} 
\label{eq:lipschitz-theta}
\|\A_i^*(f_{-i}, \theta_i') - \A_i^*(f_{-i}, \theta_i'') \|
\leq
\frac{\kappa_i}{\sigma_i } \|\theta_i' - \theta_i''\|, \; \forall \theta_i', \theta_i'' \in \Theta_i.
\end{equation}
\end{theorem}

\noindent
\textbf{Proof.}
Similar to the proof in Theorem~\ref{thm:lipschitz-f}, under Assumptions~\ref{assu:u_i} - \ref{ass:equi}(a),
$-\nabla_{a_i} \phi_i\left(a_i, f_{-i}, \theta_i\right)$ 
is strongly monotone in $a_i$ uniformly with respect to 
$f_{-i}, \theta_i$.
Then for each fixed $f_{-i}\in\F_{-i}$ and $\theta_i\in\Theta_i$, $\A_i^*(f_{-i},\theta_i)$ satisfies the variational inequality in \eqref{eq:phi-vi}.
Therefore, by fixing up $f_{-i}$ and treating 
$\theta_i$ as a varying parameter, 
Lemma~\ref{Lem:Lip-slv-mon} applies to the strongly monotone operator
$-\nabla_{a_i} \phi_i\left(\cdot, f_{-i}, \theta_i\right)$ over $\A_i$ and we obtain
\begin{align}
\nonumber
\|\A_i^*(f_{-i}, \theta_i') - \A_i^*(f_{-i}, \theta_i'') \|
& \leq 
\frac{1}{\sigma_i }
\|
\nabla_{a_i} \phi_i\left(\A_i^*(f_{-i}, \theta_i''), f_{-i}, \theta_i''\right)
-
\nabla_{a_i} \phi_i(\A_i^*(f_{-i}, \theta_i''), f_{-i}, \theta_i')
\| \\
& \leq
\frac{\kappa_i}{\sigma_i } \|\theta_i' - \theta_i''\|,
\end{align}
where the second inequality is due to the Lipschitz continuity (see Proposition~\ref{prop:lipschitz}).
$\hfill \square$

The main condition that we impose is the uniform strong concavity condition of  
$u_i(a_i, a_{-i}, \theta)$ in $a_i$. 
Guo et al.~\cite{guo2021existence} use the condition to derive
equicontinuity of the optimal response function. 
Since their derivation is based on the second-order growth condition, the optimal response function is shown to be $\frac{1}{2}$-H\"older continuous in $\theta_i$.
By contrast, our approach is based on the implicit function theorem applied to the first-order optimality condition.
This allows us to exploit the same strong concavity assumption in a different way and obtain the stronger estimate in \eqref{eq:lipschitz-theta}, that is, Lipschitz continuity of the optimal response with respect to $\theta_i$.

\section{Existence and uniqueness}
\label{se:Existence and uniqueness}

With the preparations in the preceding section, we are now ready to discuss existence and uniqueness of BNEs.

\subsection{Measurable BNE and continuous BNE}

We begin with the one that is not necessarily continuous.
\begin{theorem}[Existence and uniqueness of a measurable BNE]
\label{thm:BNE-contraction}
Consider the BNE model~\eqref{eq:BNE}. Suppose that Assumptions~\ref{assu:u_i} - \ref{ass:equi} hold.
Let $M=[M_{ij}]_{i,j=1}^n$ be the nonnegative matrix defined by
    \begin{equation}
    \label{eq:M_matrix}
    M_{ij}:=
    \begin{cases}
    0, & i=j,\\[2pt]
    \tau_{ij}/\sigma_i, & i\neq j.
    \end{cases}
    \end{equation}
Assume that there exists $w =(w_1,\dots,w_n)\in \R_{++}^n$ such that
$\|D^{-1}MD\|_\infty < 1$,
where 
$D:=\mathrm{diag}(w)$ and
$\| \cdot \|_\infty$ 
denotes the matrix infinity norm.
Then the following assertions hold.

\begin{itemize}
\item[(i)] For every $p\in[1,\infty]$, the operator $\Psi$ is a contraction on
$(\mathcal F,\|\cdot\|_{\infty,w,\mathcal L^p(\eta)})$ and
there exists a  
BNE 
which is unique 
$\eta_i$-a.e. 
for $i\in N$.

\item[(ii)]  The operator $\Psi$ is also a contraction on
$(\mathcal F,\|\cdot\|_{\infty,w})$
and the equilibrium is unique pointwise on $\Theta$.

\item[(iii)] If Assumption~\ref{assu:u_i}(c) is strengthened by requiring that the conditional distribution
$\eta_i(\cdot\mid\theta_i)$ is \underline{setwise continuous} in $\theta_i$ (i.e., for every Borel set $S \in \mathscr{B}(\Theta_{-i})$, $\eta_{i}(S \mid \theta_i)$ is continuous in $\theta_i$), then the unique equilibrium is continuous, and hence the model admits a unique CBNE.

\end{itemize}

\end{theorem}

It might be helpful to give some comments on the results before providing a proof. 
The existence and uniqueness established in
Part~(i) covers a wider range of equilibria where
each player's optimal response function is unique
over a subset of its domain with probability $1$,
allowing deviations over any subset with zero probability.
By contrast, Part~(ii) requires the optimal response function to be unique at each scenario, excluding existence of any other equilibrium. As such, we can view the equilibrium as one of the equilibria in Part (i).
From practical application perspective, 
the latter is stronger and more desirable in the sense that
a unique equilibrium is secured in any scenarios rather than in almost every scenario. 
It is important to note that under the conditions of the theorem, there is only one pointwise unique equilibrium. The almost-sure uniqueness represents
 a different angle to analyse the equilibrium, which ignores the performances of the response strategy over any set of scenarios with zero probability.

\noindent
\textbf{Proof.}
Let $\F$ be defined as in \eqref{eq:F} and  $f \in \F$. Define the optimal response mapping
\begin{equation}
\label{eq:PSI}
\Psi(f) := 
\big(
\tilde{f}_1,\dots,\tilde{f}_n
\big) 
=
\big(
\A_1^*(f_{-1}, \cdot),
\dots,
\A_n^*(f_{-n}, \cdot)
\big) \in \F.
\end{equation}
\underline{Part (i)}.
For the given vector $w:=(w_1,\dots,w_n)\in\R^n_{++}$, recall that the weighted $\mathcal{L}^p(\eta)$-norm is
$$
\|f\|_{\infty,w,\mathcal{L}^p(\eta)}
=\max_{i\in N} \left(\frac{1}{w_i}\|f_i\|_{\mathcal L^p(\eta_i)}\right),
$$
then for all $f,g\in\F$ and $p \in [1, +\infty]$, it follows that
\begin{align}
\nonumber
\|\Psi(f)-\Psi(g)\|_{\infty,w,\mathcal{L}^p(\eta)}
& =
\max_{i\in N}
\left(
\frac{1}{w_i}
\|\tilde f_i-\tilde g_i\|_{\mathcal{L}^p(\eta_i)}
\right) 
\leq
\max_{i\in N}
\left(
\frac{1}{w_i}
\sum_{j=1}^n M_{ij}\|f_j-g_j\|_{\mathcal{L}^p(\eta_j)}
\right) \\ \nonumber
& \leq
\max_{i\in N}
\left(
\frac{1}{w_i}
\sum_{j=1}^n M_{ij}w_j\,\|f-g\|_{\infty,w,\mathcal{L}^p(\eta)}
\right) \\ \nonumber
& =
\left(
\max_{i\in N}\frac{1}{w_i}
\sum_{j=1}^n M_{ij}w_j
\right)
\|f-g\|_{\infty,w,\mathcal{L}^p(\eta)}
\\
& =
\|D^{-1}MD\|_\infty
\|f-g\|_{\infty,w,\mathcal{L}^p(\eta)},
\end{align}
where the first and second inequalities follow from
\eqref{eq:lipschitz-f} and  
$\frac{1}{w_j}\|f_j - g_j\|_{\mathcal{L}^p(\eta_j)} 
\leq \|f - g\|_{\infty,w, \mathcal{L}^p(\eta)}$, respectively.
Therefore, $\Psi:\F\to\F$ is a contraction mapping with contraction modulus
\begin{equation}
\label{eq:alpha}
\alpha :=
\left(\max_{i\in N}\frac{1}{w_i}
\sum_{j =1}^n M_{ij} w_j \right)
=  
\|
D^{-1}MD 
\|_\infty
< 1.
\end{equation}
Since the product space
$\prod_{i\in N}{\cal L}^p(\Theta_i,\eta_i)$
equipped with the norm $\|\cdot\|_{\infty,w,\mathcal L^p(\eta)}$ is complete, by Banach's fixed-point theorem~\cite{zeidler1993nonlinear}, a unique BNE exists in this metric sense.
Here, uniqueness is understood with respect to this metric, that is,
two strategy profiles are identified whenever
$\|f-g\|_{\infty,w,\mathcal L^p(\eta)}=0$, or equivalently, whenever
$\|f_i-g_i\|_{\mathcal L^p(\eta_i)}=0$ for all $i\in N$.
Thus, the BNE is unique up to $\eta_i$-a.e. equality for each $i\in N$

\underline{Part (ii)}.
Recall that the weighted infinity norm is
$
\|f\|_{\infty,w} = \max_{i\in N}\frac1{w_i}\|f_i\|_\infty.
$
By the infinity norm Lipschitz estimate in \eqref{eq:lipschitz-f2},
similar to the proof in Part~(i), for any $f,g\in\mathcal F$,
\begin{align*}
\|\Psi(f)-\Psi(g)\|_{\infty,w}
=
\max_{i\in N}\frac1{w_i}\|\tilde f_i-\tilde g_i\|_\infty 
\leq
\|D^{-1}MD\|_\infty\,\|f-g\|_{\infty,w}.
\end{align*}
Therefore, $\Psi$ is also a contraction on $(\mathcal F,\|\cdot\|_{\infty,w})$.
Since $(\mathcal F,\|\cdot\|_{\infty,w})$ is complete, Banach's fixed-point theorem
yields a unique fixed point in $\mathcal F$.
Moreover, $\|f-g\|_{\infty,w}=0$ if and only if $\|f_i-g_i\|_\infty=0$ for all $i\in N$,
hence the equilibrium is unique pointwise on $\Theta$.

\underline{Part (iii)}.
If the conditional distribution
$\eta_i(\cdot\mid\theta_i)$ is setwise continuous in $\theta_i$, then for each $f_{-i}\in\F_{-i}$, 
$\phi_i(a_i,f_{-i},\theta_i)$ is jointly continuous in $(a_i,\theta_i)$, see Su and Xu~\cite[Proposition 3.1]{SU2025BNRE}.
By Berge's maximum theorem~\cite{berge1877topological,aliprantis2006infinite}, the optimal response function
$\A_i^*(f_{-i},\cdot)$ is continuous on $\Theta_i$.
Therefore, $\Psi(\mathcal F)\subset \mathcal C$.
Since the unique equilibrium $f^*$ satisfies $f^*=\Psi(f^*)$, it follows that $f^*\in\mathcal C$.
Hence the model admits a unique CBNE.
$\hfill \square$

From the proofs of Theorems~\ref{thm:lipschitz-f} and \ref{thm:BNE-contraction}, we know that the contraction property of $\Psi$ does not rely on the Lipschitz continuity of $\nabla_{a_i} \phi_i(a_i, f_{-i}, \theta_i)$ and 
$\A_i^*(f_{-i}, \theta_i)$ w.r.t. $\theta_i$.
What is essential for the contraction argument is the Lipschitz continuity of the optimal-response mapping with respect to the rivals' strategy profile.
Hence, even if the optimal responses are not continuous in types, one may still work on the measurable strategy space $\mathcal F$ and apply Banach's fixed-point theorem to obtain existence and uniqueness of a (not necessarily continuous) Bayesian Nash equilibrium in the metric induced by $\|\cdot\|_{\infty,w,\mathcal L^p}$ and $\|\cdot\|_{\infty,w}$.
The existence and uniqueness results differ from 
\cite{vives1990nash,athey2001single,mcadams2003isotone,vives2007monotone,mason2007existence,reny2011existence} where existence relies heavily on the monotonicity of the optimal response function (there is no uniqueness result in their settings with exception \cite{mason2007existence}). 
They also differ from 
\cite{meirowitz2003existence,ui2016bayesian,guo2021existence,tao2024generalized,SU2025BNRE} where the existence and uniqueness are established only for
CBNE, and the uniqueness is derived under the ex-ante 
framework. Indeed, the CBNE result established 
in Theorem~\ref{thm:BNE-contraction}(iii)
does not require equicontinuity of the optimal response functions, a key condition specified in \cite{meirowitz2003existence} and the follow-up works.

Finally, we note that the key 
contraction condition $\|D^{-1}MD\|_\infty < 1$ in \eqref{eq:alpha} may be presented in other forms. The next proposition states this.

\begin{proposition}[Equivalent contraction conditions]
\label{prop:rhoM<1}
Let $M\in \R^{n\times n}$ be a nonnegative matrix (entrywise). Then the following statements are equivalent:
\begin{enumerate}

\item[(i)] There exists a vector $w \in \R^n_{++}$ such that, with $D:=\mathrm{diag}(w)$,
$\|D^{-1}MD\|_\infty<1$;

\item[(ii)] The spectral radius of $M$ satisfies $\rho(M)<1$;

\item[(iii)] There exists a vector $w\in \R^n_{++}$ such that
$Mw<w$ \text{(componentwise)}.

\end{enumerate}
Moreover, if $\rho(M)<1$, then one may choose
$w:=\sum_{k=0}^{\infty} M^k\mathbf{1},$
for which $w\in \R^n_{++}$ and $Mw  < w$.
\end{proposition}

\noindent
\textbf{Proof.}
\underline{(i)$\Rightarrow$(ii)}. 
Since $D=\mathrm{diag}(w)$ is non-singular, then
the matrix $D^{-1}MD$ is similar to $M$, which implies that $\rho(D^{-1}MD) = \rho(M)$.
By the standard bound 
$\rho(A) \leq 
\|
A
\|$ 
for any induced matrix norm 
$\|
\cdot
\|$
(see Horn and Johnson~\cite[Theorem 5.6.9]{horn2012matrix}), we obtain
$$
\rho(M) = \rho(D^{-1}MD) \leq 
\|
D^{-1}MD
\|_\infty < 1.
$$

\underline{(ii)$\Rightarrow$(iii)}. 
Suppose $\rho(M)<1$. Then the Neumann series $\sum_{k=0}^{\infty} M^k$ converges to $(I - M)^{-1}$.
Define 
$$
w:= (I-M)^{-1}\mathbf{1} = \sum_{k=0}^{\infty} M^k\mathbf{1}.
$$
Since $M\ge 0$, for each $k \in \mathbb{N}$, the vector $M^k\mathbf 1$ is nonnegative. 
Therefore,  
$w = \mathbf{1} + \sum_{k=1}^{\infty} M^k\mathbf{1} \geq \mathbf{1}$, which implies that $w \in \R_{++}^n$.
Moreover,
$(I-M)w=\mathbf{1}$, and hence $Mw = w-\mathbf{1} < w$.

\underline{(iii)$\Leftrightarrow$(i)}. 
Note that $(D^{-1}MD)_{ij}=M_{ij}w_j/w_i$. Then it follows that
$$
\|
D^{-1}MD
\|_\infty
=\max_{i\in N}\sum_{j=1}^n 
\left \vert (D^{-1}MD)_{ij} \right \vert
=\max_{i\in N}\frac{1}{w_i}\sum_{j=1}^n M_{ij}w_j
=\max_{i\in N}\frac{(Mw)_i}{w_i},
$$
where the second equality is due to the fact that $M\geq 0$ and $w \in \R_{++}^n$.
Therefore,
$\|
D^{-1}MD
\|_\infty<1$
if and only if $(Mw)_i<w_i$ for all $i \in N$, that is, $Mw < w$.
$\hfill\square$


\begin{remark} 

The condition $\rho(M)<1$ also admits a natural interpretation in terms of strategic feedback. 
Indeed, for $i\neq j$, the quantity $M_{ij}=\tau_{ij}/\sigma_i$ compares two marginal-utility effects: 
$\tau_{ij}$ is an upper bound on the effect of player $j$'s action on player $i$'s marginal utility, whereas $\sigma_i$ is a lower bound on the sensitivity of player $i$'s marginal utility with respect to its own action. 
Thus, each entry of $M$ represents the ratio between a cross-player marginal effect and player $i$'s own marginal sensitivity, and the matrix $M$ captures these interaction bounds across players.

From a dynamic viewpoint, the condition $\rho(M)<1$ ensures that there exists $w\in \R_{++}^n$ such that, with $D=\mathrm{diag}(w)$, $\alpha=\|D^{-1}MD\|_\infty<1.$
Consequently, after $k$ successive rounds of optimal-response updates, the aggregate indirect feedback is bounded by $\alpha^k$, and hence converges to zero geometrically as $k\to\infty$. 
In economic terms, indirect strategic effects become progressively weaker rather than amplifying over time. 
This is precisely the mechanism underlying the contraction property of the optimal-response mapping.

Related interpretations based on interaction matrices and repeated feedback also appear in the network games literature, see, e.g., \cite{ballester2006s,parise2019variational}. 
This viewpoint is also consistent with the broader economic networks framework of \cite{acemoglu2015networks,acemoglu2015systemic}, where weighted directed network structures are used to represent interconnections and to characterize the propagation of direct and indirect effects through the network.

\end{remark}

Note that the contraction condition presented here
differs from the uniqueness condition used by Ui~\cite{ui2016bayesian} and the follow-up works, the next example illustrates this.

\begin{example}
\label{ex:asym_2p_BG}
Consider an asymmetric two-player Bayesian game,
where
$$
u_1(a_1,a_2,\theta)
=
-\theta_1 a_1-\frac12 a_1^2-2a_1a_2,
\quad
u_2(a_1,a_2,\theta)
=
-\theta_2 a_2-\frac12 a_2^2-\frac25 a_1a_2,
$$
and $\A_i=[0,\overline a],\
\Theta_i=[\underline\theta,\overline\theta]$
for  $i=1,2$.
Observe that 
$$
\nabla_{a_1}u_1(a_1,a_2,\theta)
=
-\theta_1 - a_1-2a_2, \quad
\nabla_{a_2}u_2(a_1,a_2,\theta)
=
-\theta_2 - a_2-\frac25 a_1,
$$
which implies that $\nabla_{a_1}u_1$ is Lipschitz continuous in $a_2$ with modulus
$\tau_{12}=2$, while $\nabla_{a_2}u_2$ is Lipschitz continuous in $a_1$ with modulus
$\tau_{21}=\frac25$.
Assume that 
the type profile $\theta=(\theta_1,\theta_2)$ is drawn from a distribution $\eta$
supported on $\Theta_1\times\Theta_2$ and
the conditional distribution $\eta_i(\cdot \mid \theta_i)$ is setwise continuous in $\theta_i$.
The expected utility of player $1$ is
$$
\phi_1(a_1,f_2,\theta_1)
=
-\theta_1a_1-\frac12 a_1^2
- 2a_1
\int_{\Theta_2} f_2(\theta_2)d\eta_1(\theta_2 | \theta_1),
$$
and the expected utility of player $2$ is
$$
\phi_2(a_2,f_1,\theta_2)
=
-\theta_2 a_2-\frac12 a_2^2
-
\frac25 a_2
\int_{\Theta_1} f_1(\theta_1) d\eta_2(\theta_1 | \theta_2).
$$
Since
$\nabla_{a_1}^2\phi_1(a_1,f_2,\theta_1) =
\nabla_{a_2}^2\phi_2(a_2,f_1,\theta_2)=-1$,
then $\phi_i$, $i=1,2$, is strongly concave in $a_i$
with concavity moduli $\sigma_1=\sigma_2=1$.
Consequently, we can figure out
$
M=
\begin{pmatrix}
0 & 2 \\
\frac{2}{5} & 0
\end{pmatrix},
$
with 
$\rho(M)=\sqrt{\frac45}<1$.
By Theorem~\ref{thm:BNE-contraction}(iii), 
the Bayesian game admits a unique continuous equilibrium.

Next, we show that the strict monotonicity condition in Guo et al.~\cite[Theorem 3.5]{guo2021existence} fails for this example.
Define the payoff-gradient mapping
$$
G(a,\theta)
:=
\begin{pmatrix}
\nabla_{a_1}u_1(a_1,a_2,\theta)\\
\nabla_{a_2}u_2(a_1,a_2,\theta)
\end{pmatrix}
=
\begin{pmatrix}
-\theta_1-a_1-2a_2\\
-\theta_2-a_2-\frac25 a_1
\end{pmatrix}.
$$
For any two continuous strategy profiles $f,\tilde{f} \in \C$, let
$\Delta_i(\theta_i):=f_i(\theta_i)-\tilde f_i(\theta_i), i=1,2.$
Then it follows that
$$
G(f(\theta),\theta)-G(\tilde f(\theta),\theta)
=
\begin{pmatrix}
-\Delta_1(\theta_1)-2\Delta_2(\theta_2)\\
-\frac25\Delta_1(\theta_1)-\Delta_2(\theta_2)
\end{pmatrix},
$$
and hence
$$
\big(G(f(\theta),\theta)-G(\tilde f(\theta),\theta)\big)^\top
\big(f(\theta)-\tilde f(\theta)\big)
=
-\Delta_1(\theta_1)^2-\Delta_2(\theta_2)^2-\frac{12}{5}\Delta_1(\theta_1)\Delta_2(\theta_2).
$$
Now we choose 
$f_1(\theta_1)\equiv \epsilon,
f_2(\theta_2)\equiv 0,$
and
$\tilde{f}_1(\theta_1)\equiv 0,\
\tilde{f}_2(\theta_2)\equiv \epsilon,$
where $0<\epsilon<\overline a$.
Then we have
$\Delta_1(\theta_1)\equiv \epsilon,\
\Delta_2(\theta_2)\equiv -\epsilon$,
and it follows that
$$
\big(G(f(\theta),\theta)-G(\tilde f(\theta),\theta)\big)^\top
\big(f(\theta)-\tilde f(\theta)\big)
=
\frac25 \epsilon^2>0,
\quad
\forall \theta\in\Theta.
$$
Therefore,
$$
\int_\Theta
\big(G(f(\theta),\theta)-G(\tilde f(\theta),\theta)\big)^\top
\big(f(\theta)-\tilde f(\theta)\big)\,\eta(d\theta)
=
\frac25 \epsilon^2>0,
$$
which violates the strict monotonicity condition in Guo et al.~\cite[Theorem 3.5]{guo2021existence}.
Hence, the uniqueness result in Guo et al.~\cite{guo2021existence} does not apply to this example.
Moreover, Guo et al.~\cite{guo2021existence} explicitly note that its uniqueness condition is weaker than the Jacobian-based negative-definiteness condition in Ui~\cite[Lemma 4 and Proposition 2]{ui2016bayesian}.
Therefore, the latter condition is not satisfied either.
Indeed, the Jacobian matrix of $G(a,\theta)$ with respect to $a=(a_1,a_2)$ is
$J(a,\theta)
=
\begin{pmatrix}
-1 & -2 \\
-\frac25 & -1
\end{pmatrix}$,
and its symmetric part is
$\frac{J(a,\theta)+J(a,\theta)^\top}{2}
=
\begin{pmatrix}
-1 & -\frac65 \\
-\frac65 & -1
\end{pmatrix}$.
Since 
$\det 
\begin{pmatrix}
-1 & -\frac65 \\
-\frac65 & -1
\end{pmatrix}
= 1-\Big(\frac65\Big)^2<0$,
it shows that the matrix $\frac{J(a,\theta)+J(a,\theta)^\top}{2}$ is not negative definite.
Therefore, the sufficient condition in Ui~\cite[Lemma 4 and Proposition 2]{ui2016bayesian} also fails to ensure the uniqueness.
\end{example}

\subsection{Lipschitz continuous BNE}

Theorem~\ref{thm:BNE-contraction} establishes existence and uniqueness of a BNE and, under an additional continuity requirement on the conditional distributions, also yields existence and uniqueness of a CBNE.
However, these conclusions are still qualitative in nature: they do not quantify how the equilibrium strategy varies with the player's own type parameter $\theta_i$.
To obtain such quantitative regularity, we impose additional conditions ensuring Lipschitz dependence on $\theta_i$.
The next theorem strengthens the preceding CBNE result by establishing existence of a Lipschitz CBNE and, under the contraction condition, its uniqueness.
\begin{theorem}[Existence and uniqueness of Lipschitz CBNE]
\label{thm:exist-unique-cbne}
Consider the BNE model~\eqref{eq:BNE}.
Suppose that Assumptions \ref{assu:u_i} - \ref{ass:lipschitz} hold.
Then the following assertions hold.

\begin{itemize}
    \item [(i)] The problem has a CBNE with all players' optimal 
    response functions being uniformly Lipschitz continuous at the equilibrium.

    \item [(ii)] 
    If, in addition, the cross-effect matrix $M$ defined in \eqref{eq:M_matrix}
    satisfies $\rho(M)<1$, then the CBNE is unique.

\end{itemize}
 
\end{theorem}

\noindent
\textbf{Proof.}
Let $\F$ and $\Psi(\cdot)$ be defined as in \eqref{eq:F} and \eqref{eq:PSI}, respectively.

\underline{Part (i)}.
This result 
follows the classical Schauder’s fixed-point theorem (cf. Lemma~\ref{fixedpoint}) route used in~\cite{guo2021existence}, where the equicontinuity of each player's response function is shown to be 
$\frac{1}{2}$-H\"older continuous (when the order of growth is $2$). Here we use the Lipschitz continuity 
of the response functions as 
stated in \eqref{eq:lipschitz-theta}.
Let $\C$ be defined as in (\ref{eq:F-C})
and equip it with the topology induced by the infinity norm $\| \cdot \|_\infty$.
We verify that $\Psi: \C \rightarrow \C$ is a compact operator.
By Theorem~\ref{thm:lipschitz-theta},
$\A_i^*(f_{-i}, \theta_i)$ 
is equicontinuous over
$\Theta_i$. 
Since $\bigcup_{\theta_i \in \Theta_i} \A_i^*(f_{-i}, \theta_i) \subset 
\A_i$ is uniformly bounded, 
then by the Arzel\`a–Ascoli theorem (cf. Lemma \ref{ascoli}),
$\Psi(\C)$ is relatively compact. 
Moreover, 
it follows by  Theorem~\ref{thm:lipschitz-f} that 
$\C_{-i} \rightarrow \C_i, f_{-i} \mapsto \A^*_i(f_{-i}, \cdot)$ is continuous in $f_{-i}$
for $i \in N$. 
This shows that $\Psi: \C \rightarrow \C$ is a continuous operator. 
By Lemma \ref{fixedpoint}, 
problem (\ref{eq:BNE}) has 
a CBNE.
Let $f^*$ be a CBNE.
For any $\theta', \theta'' \in \Theta$, 
\begin{align*}
\| f^*(\theta') - f^*(\theta'')\| 
& = \sqrt{\sum_{i =1}^n 
\left \|\A_i^*(f^*_{-i}, \theta_i') - \A_i^*(f^*_{-i}, \theta_i'') \right \|^2 } \\
& \leq 
\sqrt{ \sum_{i =1}^n
\left( \frac{\kappa_i}{\sigma_i} \right)^2
\|\theta_i' - \theta_i''\|^2} \quad (\inmat{by \eqref{eq:lipschitz-theta}})
\\
&\leq
\left( \max_{i \in N} \frac{\kappa_i}{\sigma_i} \right)
\|\theta' - \theta''\|,
\end{align*}
which implies that
$f^*$ is Lipschitz continuous in $\theta$ with modulus 
$ \max\limits_{i \in N} \frac{\kappa_i}{\sigma_i}$.

\underline{Part (ii)}.
Since Assumption~\ref{ass:lipschitz}(c) implies that
$\eta_i(\cdot\mid\theta_i)$ is setwise continuous in $\theta_i$ for each $i\in N$, by Theorem~\ref{thm:BNE-contraction}(iii) and Proposition~\ref{prop:rhoM<1}, the Lipschitz CBNE is unique pointwise on $\Theta$. 

$\hfill \square$

The following symmetric three-player Bayesian Cournot example illustrates that our 
contraction condition $\rho(M) < 1$ can naturally hold in practice.

\begin{example}[Symmetric three-player Bayesian Cournot game]
\label{example:3p_sym_BCG}
Consider a symmetric Bayesian Cournot game with three players.  
For each player $i=1,2,3$, let
$$
u_i(a_i,a_{-i},\theta) 
= a_i \, p(a_1,a_2, a_3) - c_i(a_i,\theta_i),
$$
where $p(a_1,a_2,a_3) = \alpha - \beta(a_1+a_2+a_3), \ \alpha,\beta>0,$ is the inverse demand function
and $c_i(a_i,\theta_i) = \theta_i a_i + \frac{c}{2}a_i^2,\ c>0$, is the cost function.
The action and type spaces are
$\A_i = [0,\overline{a}], \ 
\Theta_i = [\underline{\theta},\overline{\theta}], i=1,2, 3,$
 and the type profile $\theta=(\theta_1,\theta_2, \theta_3)$ is drawn from a distribution $\eta$ supported on $\Theta_1\times \Theta_2 \times \Theta_3$. Assume the density function $q_i(\theta_{-i}|\theta_i)$ of the conditional distribution $\eta_i(\theta_{-i}|\theta_i)$ is Lipschitz continuous in $\theta_i$ with modulus $\gamma_i>0$.
Then the utility function of player $i$ takes the form
$$
u_i(a_i,a_{-i},\theta)
= (\alpha-\theta_i)a_i - \Big(\beta+\frac{c}{2}\Big)a_i^2 - \sum_{j \neq i} \beta a_i a_j.
$$
The gradient of $u_i$ is
$$
\nabla_{a_i} u_i(a_i,a_{-i},\theta)
= (\alpha-\theta_i) - (2\beta+c)a_i - \sum_{j \neq i} \beta a_j,
$$
which is Lipschitz continuous in $a_{j}$ with modulus $\tau_{ij}=\beta,\ j \neq i$.
And the expected utility of player $i$ given rival strategy $f_{-i}$ is
$$
\phi_i(a_i,f_{-i},\theta_i) 
= (\alpha-\theta_i)a_i - \Big(\beta+\frac{c}{2}\Big)a_i^2 
-  \beta a_i \int_{\Theta_{-i}} 
\left( \mathbf{1}_2^\top
f_{-i}(\theta_{-i}) \right)
q_i(\theta_{-i}|\theta_i)\, d\theta_{-i},
$$
where $\mathbf{1}_2 := (1, 1)^\top$.
It is easy to verify that $\phi_i$ is strongly concave in $a_i$, since 
$$
\nabla^2_{a_i}\phi_i(a_i,f_{-i},\theta_i) = -(2\beta+c),
$$
which gives a strong concavity modulus $\sigma_i=2\beta+c$. Furthermore, the gradient
$$
\nabla_{a_i}\phi_i(a_i,f_{-i},\theta_i) 
= \alpha - \theta_i 
- \beta \int_{\Theta_{-i}} 
\left( \mathbf{1}_2^\top
f_{-i}(\theta_{-i}) \right)
q_i(\theta_{-i}|\theta_i)\, d\theta_{-i}
- (2\beta+c)a_i
$$
is Lipschitz continuous in $\theta_i$ with modulus  
$$
\kappa_i = 1 + \beta\gamma_i \int_{\Theta_{-i}} \mathbf{1}_2^\top
f_{-i}(\theta_{-i})\, d\theta_{-i} 
\leq 1 + 2\beta\gamma_i \overline{a}.
$$
Since $\tau_{ij}=\beta$ and $\sigma_i=2\beta+c$, 
we have
$$
\|M\|_\infty =
\max_{i\in N}\sum_{j\neq i}\frac{\tau_{ij}}{\sigma_i}
= \frac{2\beta}{2\beta+c} <1.
$$
Since
$\rho(M)\leq \|M\|_\infty <1,$
by Theorem~\ref{thm:exist-unique-cbne}, it guarantees that the optimal response mapping is a contraction mapping. Therefore, the Bayesian Cournot game admits a unique equilibrium in which all players’ strategy functions are uniformly Lipschitz continuous.
\end{example}

Before concluding the section, we would like to add a note that we might use the contraction property of $\Psi$ defined in \eqref{eq:PSI} to calculate an approximate BNE via an iterative process 
\begin{equation}
\label{eq:contract-Phi-f}
f^{k+1} =\Psi(f^k), \quad \text{for }
k \in \mathbb{N}.
\end{equation}
The scheme requires solving $N$ optimization problems at each iteration. Since $\|f^k-f^*\|_{\infty,w} \leq \alpha^k\|f^0-f^*\|_{\infty,w}$, then we can estimate the number of iterations to reach a prescribed precision.  
Unlike the existing computational methods based on variational inequality (see e.g.,~\cite{guo2021existence,tao2024generalized,SU2025BNRE}), 
this approach does not require $f$ to be continuous. Of course, we need to discretize the set $\Theta$ in advance.
In the case when $f$ is Lipschitz continuous, we may use piecewise-linear approximation, i.e., by confining $f$ to piecewise-linear functions defined over a specified set of breakpoints in the iterative scheme \eqref{eq:contract-Phi-f}. As such, we may obtain in the end a piecewise-linear estimate of $f^*$, and we can easily derive an error bound for the estimate. We omit the details as this is not the focus of this paper.

\section{BNE under monotonicity}
\label{se:Monotone}

In Section~\ref{se:Lipschitz}, we 
discuss the case where each player's optimal response function is Lipschitz continuous 
with respect to variations of their types and rivals' actions. 
In some practical applications, 
the optimal response function $\A_i^*(f_{-i},\theta_i)$
may display some monotonicity w.r.t.~variation type parameter $\theta_i$ and the rivals' actions.
For example, in Bertrand competition, a higher marginal cost (type) pushes up the optimal price (action), and rivals’ higher prices also raise a firm’s optimal response price (strategic complements), see~\cite{vives1990nash, anderson1992discrete}.
This kind of research is in alignment with the literature of monotone comparative statics, 
see e.g.,~\cite{milgrom1994monotone,  topkis1998supermodularity,athey2002monotone}, which uses order-theoretic tools such as quasi-supermodularity and single-crossing property to characterize how optimal decisions react to changes in parameters.
We briefly recall some relevant notions.

A \emph{partially ordered set} (poset) is a pair $(S,\geq)$, where $\geq$ is a reflexive, transitive, and antisymmetric binary relation on $S$. For $x,y\in S$, we write $x>y$ if $x\geq y$ and $x\neq y$.
A poset $(S,\geq)$ is (totally) ordered if for any $x,y\in S$, either $x\geq y$ or $y\geq x$.
Given a non-empty subset $L\subset S$, an element $u \in S$ is an upper bound of $L$ if $u \geq a$ for all $a\in L$. If there is an upper bound $u^*$ of $L$ such that $u^* \leq u$ for every upper bound $u$ of $L$, then $u^*$ is the \emph{least upper bound} (or supremum) of $L$, denoted $\vee L$. The \emph{greatest lower bound} (or infimum) $\wedge L$ is defined analogously. 

A poset $(S,\geq)$ is a \emph{lattice} if any two elements
$x,y\in S$ have a least upper bound $x\vee y:=\vee\{x,y\}$ and a greatest
lower bound $x\wedge y:=\wedge\{x,y\}$ in $S$. 
A lattice is \emph{complete}
if every non-empty subset $L\subset S$ has a supremum $\vee L$ and an infimum $\wedge L$ in $S$.
For example, $S = \{(0,0), (0,1), (1,0)\} \subset \R^2$ is not a lattice, since $(0,1)$ and $(1,0)$ have no joint upper bound in $S$.
A subset $L$ of the lattice $S$ is a \emph{sublattice} of $S$ if the supremum and infimum of any two elements of $L$ also belong to $L$.
If $S$ is a complete lattice, a sublattice $L\subset S$ is called \emph{subcomplete} if, for every non-empty subset $E\subset L$, both $\vee E$ and $\wedge E$, computed in $S$, belong to $L$.

Given a lattice $(S,\geq,\vee,\wedge)$ and an
index set $\Theta$, a function $f:S\times\Theta\to\R$ 
is \emph{quasi-supermodular} 
in $x$ if
$$
f(x' , \theta) \geq (>) f(x'\wedge x ,\theta)
\Rightarrow
f(x'\vee x ,\theta) \geq (>) f(x ,\theta)
$$
for all $x', x \in S$ and all $\theta \in \Theta$
(weak inequality implies weak inequality and strict
inequality implies strict inequality).
The function $f$ is \emph{supermodular} in $x$
if 
$$
f(x'\vee x ,\theta)
+
f(x'\wedge x ,\theta)
\geq
f(x' , \theta) + f(x ,\theta),
\ \forall x', x \in S, \theta \in \Theta.
$$
If $f$ is positive, it is 
\emph{log-supermodular} in $x$
if 
$$
f(x'\vee x ,\theta)
\cdot
f(x'\wedge x ,\theta)
\geq
f(x' , \theta) 
\cdot
f(x ,\theta),
\ \forall x', x \in S, \theta \in \Theta.
$$
Both properties are stronger than quasi-supermodularity.
In the case that $S = \R^m$ and $f$ is twice-continuously differentiable in $x$, then $f$ is supermodular if and only if $\frac{\partial^2 f(x,\theta)}{\partial x_i\,\partial x_j}\geq 0$
for all $x \in S, \theta \in \Theta$ and $i \neq j$, 
see~\cite[Theorems 3.1 and 3.2]{topkis1978minimizing}.

Given a parameter poset $(T,\geq)$, a function
$g: S\times T\times\Theta\to\R$ is said to satisfy \emph{single-crossing property}
in $(x,t)$ if
$$
g(x', t , \theta) \geq (>) g(x, t, \theta)
\Rightarrow
g(x', t' , \theta) \geq (>) g(x, t', \theta)
$$
for all $x' \geq x \in S$,
all $t' \geq t \in T$, and all
$\theta \in \Theta$.
A sufficient condition of the single-crossing property is 
\emph{increasing differences}
in $(x,t)$, that is,
$g(x', t, \theta) - g(x, t, \theta)$ is increasing in $t$ for $x' \geq x$, or equivalently, if 
$g(x, t', \theta) - g(x, t, \theta)$ is increasing in $x$ for $t' \geq t$. 
Supermodularity is a stronger property than increasing differences: if $T$ is also a lattice and if $g$ is supermodular over $X \times T$, then $g$ has increasing differences in $(x,t)$. 
In addition, if the positive function $g$ is log-supermodular over $X \times T$, it also satisfies the weak single-crossing property in $(x,t)$.

Reny~\cite{reny2011existence} further introduces the concepts of \emph{weakly quasi-supermodular} and 
\emph{weak single-crossing}, where both of them release the requirements for strict inequalities, playing important roles in the following discussion.

\begin{assumption}
\label{ass:lattice}
For each player $i \in N$, the action space $\A_i$ is a sublattice of $\R^{z_i}$.
\end{assumption}

In Section~\ref{se:The model}, we only assume that $\A_i\subset\R^{z_i}$ is compact and convex, 
which does not ensure the lattice property under the coordinatewise (natural) order.
For instance, let $\A_i := \left\{(a,b) \in \R^2: a^2 + b^2 \leq 1 \right\}$. 
Then $\A_i$ is compact and convex, but it is not a sublattice of $\R^2$: with $x=(1,0)$ and $y=(0,1)$ in $\A_i$, the coordinatewise join is $x\vee y=(1,1) \notin \A_i$. 
Thus, we cannot start monotone iterations from a bottom/top action, nor guarantee that the pointwise infimum/supremum of monotone sequences of actions remain in $\A_i$.

By Topkis~\cite[Theorem~2.3.1]{topkis1998supermodularity}, a sublattice of $\R^{z_i}$ is subcomplete if and only if it is compact. Hence any non-empty compact sublattice is a complete lattice under the coordinatewise order.
A simple and widely used instance is the closed interval in $\R^{z_i}$, that is,
$$\mathcal{A}_i:= 
\left \{ 
a_i = (a_{i1},\dots,a_{iz_i}) 
\in \R^{z_i}
\ \big| \  
\underline{a}_{ij}
\leq 
a_{ij}
\leq 
\overline{a}_{ij},\ 
\inmat{for}\; j=1,\cdots, 
z_i 
\right \},
$$ 
Assumption~\ref{ass:lattice} strengthens the structure to a complete lattice precisely to guarantee the presence of extremal elements (so that we can initialize from $\wedge\A_i$ or $\vee\A_i$) and the closure of $\A_i$ under arbitrary supremum and infimum.  
These properties are crucial for constructing extremal equilibria via monotone optimal–response iterations and for carrying out lattice–theoretic monotone comparative statics.
In our model, we interpret vector inequalities in the componentwise sense, that is, for any $m \in \mathbb{N}$ and $x , y \in \R^{m}$,
$
x \geq y \ \Longleftrightarrow\ 
x_{j} \geq y_{j},\ \forall j=1,\dots,m.
$
The corresponding supremum and infimum are defined componentwise by
$$
(x \vee y)_j := \max\{x_{j}, y_{j} \},\quad
(x \wedge y)_j := \min\{ x_{j}, y_{j} \},\quad 
\text{for}\;j=1,\dots,m.
$$
We then endow each strategy space $\mathcal{F}_i$ with the pointwise partial order induced by the above order on $\R^{z_i}$, i.e., for any $f_i,g_i \in \mathcal{F}_i$,
$
f_i \succeq g_i \
\Longleftrightarrow\
f_i(\theta_i)\geq g_i(\theta_i),\ \forall \theta_i\in \Theta_i.
$
Under this order, $(\mathcal{F}_i,\succeq)$ is a lattice. Indeed, for any $f_i,g_i\in\mathcal{F}_i$, we define
$$
(f_i\vee g_i)(\theta_i):= f_i(\theta_i)\vee g_i(\theta_i),
\quad
(f_i\wedge g_i)(\theta_i):= f_i(\theta_i)\wedge g_i(\theta_i),
\quad \forall \theta_i\in \Theta_i.
$$
Since $f_i$ and $g_i$ are measurable and the mappings $(x,y)\mapsto x\vee y$ and $(x,y)\mapsto x\wedge y$ are measurable, both $f_i\vee g_i$ and $f_i\wedge g_i$ are measurable.
Moreover, for each $\theta_i\in\Theta_i$,
$(f_i\vee g_i)(\theta_i),\ (f_i\wedge g_i)(\theta_i)\in \mathcal{A}_i$, hence we have $f_i\vee g_i,\ f_i\wedge g_i\in\mathcal{F}_i$.
We further equip $\mathcal{F}_{-i}$ and $\mathcal{F}$ with the componentwise (product) order, i.e., for $f_{-i},g_{-i}\in\mathcal{F}_{-i}$,
$
f_{-i}\succeq g_{-i}
\ \Longleftrightarrow\
f_j \succeq g_j,
\ \forall j\neq i,
$
and similarly, for $f,g\in\mathcal{F}$, $f\succeq g$ iff $f_i\succeq g_i$ for all $i\in N$.
With these orders, $(\mathcal{F}_{-i},\succeq)$ and $(\mathcal{F},\succeq)$ are lattices.

Based on the preliminaries above, we are ready to establish monotone comparative statics results for the optimal response function.

\begin{proposition}[Monotonicity of the optimal response]
\label{prop: monotonicity-f}
Consider the BNE model \eqref{eq:BNE}.
Suppose that Assumptions~\ref{assu:u_i} - \ref{ass:lattice} hold, and for $i \in N$, $\phi_i(a_i,f_{-i},\theta_i)$ 
is weakly quasi-supermodular in $a_i \in \A_i$
for all $f_{-i} \in \F_{-i}$
and $\theta_i \in \Theta_i$. 
Then the following assertions hold.

\begin{itemize}
\item [(i)] For each fixed $\theta_i \in \Theta_i$, 
if $\phi_i(a_i, f_{-i}, \theta_i)$ is weak single-crossing over $\A_i \times \F_{-i}$,
then the optimal response function $\A_i^*(f_{-i}, \theta_i)$ is Lipschitz continuous and increasing in $f_{-i}\in \F_{-i}$, that is, for any $f_{-i} \succeq g_{-i}$, $\A_i^*(f_{-i}, \theta_i) \geq \A_i^*(g_{-i}, \theta_i)$ and $\|\A_i^*(f_{-i}, \theta_i) - \A_i^*(g_{-i}, \theta_i) \| \leq \sum_{j \neq i} \frac{\tau_{ij}}{\sigma_i}\|f_j - g_j \|_\infty$. 

\item [(ii)] For each fixed $f_{-i}\in \F_{-i}^+$, 
if $\phi_i(a_i, f_{-i}, \theta_i)$ is weak single-crossing over $\A_i \times \Theta_{i}$,
then the optimal response function $\A_i^*(f_{-i}, \theta_i)$ is Lipschitz continuous and increasing in $\theta_i\in \Theta_i$, that is, for any $\theta_i \geq \theta_i'$, $\A_i^*(f_{-i}, \theta_i) \geq \A_i^*(f_{-i}, \theta_i')$ and 
$\|\A_i^*(f_{-i}, \theta_i) - \A_i^*(f_{-i}, \theta_i') \| \leq \frac{\kappa_i}{\sigma_i}\|\theta_i - \theta_i'\|$. 
\end{itemize}
\end{proposition}

The Lipschitz estimates in parts (i) and (ii) follow from
Theorems~\ref{thm:lipschitz-f} and~\ref{thm:lipschitz-theta}, respectively.
By virtue of monotone comparative statics results as shown in \cite[Proposition~3.2]{SU2025BNRE}, the above assertions follow directly, thus we omit the proof here.
Sufficient conditions for quasi-supermodularity and single-crossing in Bayesian games are well-documented such as supermodularity and increasing differences (see e.g.,~\cite{athey2001single,vives2007monotone}),
and the log-supermodularity (e.g.,~\cite{athey2001single,athey2002monotone,john2012}).
Specifically, if $\mathcal{A}_i \subset \R$, then the supermodularity is automatically satisfied, due to the fact that $\R$ is totally ordered and hence
$(a_i'\vee a_i), (a_i'\wedge a_i)
\in \{a_i', a_i\}$ for all $a_i', a_i \in \mathcal{A}_i$.

In the next theorem, we use the monotone iteration technique of Van Zandt and Vives~\cite{vives2007monotone} to establish the existence of greatest and least equilibria. Under an additional own type single-crossing condition, these equilibria are monotone and hence become CMBNEs.

\begin{theorem}[Existence and uniqueness of
greatest and least CMBNEs]
\label{thm:exist-unique-CMBNE}
Consider the BNE model \eqref{eq:BNE}.
Suppose: 
(a) Assumptions~\ref{assu:u_i} - \ref{ass:lattice} hold; 
(b) $\phi_i(a_i, f_{-i}, \theta_i)$ 
is weakly quasi-supermodular 
over $\mathcal{A}_i$ for all $f_{-i} \in \F_{-i}$
and $\theta_i \in \Theta_i$, and is weak single-crossing over $\A_i \times \F_{-i}$ for all $\theta_i \in \Theta_i$.
Then the following assertions hold.

\begin{itemize}

    \item [(i)] 
    The problem has a greatest BNE and a least BNE, with all players' equilibrium strategy functions being uniformly Lipschitz continuous. 

    \item [(ii)] If, in addition, 
    $\phi_i(a_i, f_{-i}, \theta_i)$
    satisfies weak single-crossing over $\A_i \times \Theta_i$ for all $f_{-i} \in \F_{-i}^+$, then the greatest and the least Lipschitz CBNEs are monotone.
    
    \item [(iii)] Furthermore, 
    if
    the cross-effect matrix $M$ defined in \eqref{eq:M_matrix}
    satisfies $\rho(M)<1$,
    then the greatest and least Lipschitz CMBNEs coincide, hence the equilibrium is unique.
\end{itemize}

\end{theorem}

\noindent
\textbf{Proof.}
\underline{Part (i)}.
We only show existence of the greatest BNE 
as existence of the least BNE can be proved analogously. 
By Proposition~\ref{prop: monotonicity-f}(i), the mapping 
$\F_{-i}\ni f_{-i}\mapsto \A_i^*(f_{-i},\cdot)\in\C_i$
is increasing with respect to 
$f_{-i} \in \F_{-i}$,
that is, $\A^*_i(f'_{-i}, \cdot) \succeq \A^*_i(f''_{-i}, \cdot)
$ for all $f'_{-i} \succeq f''_{-i} \in \F_{-i}$.
Let $\Psi(f)$ be defined as in (\ref{eq:PSI}), and hence we can conclude from the discussions above that the operator $\Psi(\cdot)$ is increasing in $f \in \F$.

Next, we show existence of the greatest BNE. 
For each $i$, let $f_i^0(\theta_i):=\vee \A_i$ for all $\theta_i \in \Theta_i$, i.e., the constant function at the top element of $\A_i$. 
Let $f^0 := (f^0_1, \dots, f^0_n) \in \F^+$, and define the sequence of functions 
$$
f^k := \Psi(f^{k-1}), \quad \inmat{for} \; k \in \mathbb{N}_+.
$$ 
Since $f^0$ is the profile of the greatest strategy functions, it follows that $f^1 \preceq f^0$. 
Together with the increasing property of operator $\Psi(\cdot)$, we can deduce that
$f^2 = \Psi(f^1) \preceq \Psi(f^0) = f^1$.
By induction, it follows that the sequence of optimal response functions
$\{f^k\}_{k \in \mathbb{N}}$ is decreasing.
This implies that 
$\{ f^k_i(\theta_i)\}$ is a decreasing sequence for all $\theta_i \in \Theta_i, i \in N$.
Since the sequence is lower bounded by $\wedge\A_i$,
then every decreasing sequence in $\A_i$ converges to its infimum, denoted by $f^\infty_i(\theta_i)$.
Thus, we can conclude that $f^k_i$ 
converges pointwise to $f^\infty_i$.


By Theorem~\ref{thm:lipschitz-theta}, for each $i\in N$, the functions
$f_i^k$ are Lipschitz continuous on $\Theta_i$ with the common modulus
$\kappa_i/\sigma_i$, independent of $k$. Hence, their pointwise limit $f_i^\infty$ is also Lipschitz continuous. Indeed, for any $\theta_i,\theta_i'\in\Theta_i$, we have
$$
\|f_i^\infty(\theta_i)-f_i^\infty(\theta_i')\|
=
\lim_{k\to\infty}
\|f_i^k(\theta_i)-f_i^k(\theta_i')\|
\leq
\frac{\kappa_i}{\sigma_i}\|\theta_i-\theta_i'\|.
$$
Thus $f_i^\infty$ is continuous. Since $\Theta_i$ is compact and $\{f_i^k\}_{k\in\mathbb N}$ decreases pointwise to the continuous function $f_i^\infty$, Dini's theorem implies that
$f_i^k\to f_i^\infty$ uniformly on $\Theta_i$. Hence
$f^k\to f^\infty$ uniformly under the infinity norm.
Moreover, $\Psi(\F)\subset\C$ by Theorem~\ref{thm:lipschitz-theta}, and
$\Psi$ is continuous under $\|\cdot\|_\infty$ by Theorem~\ref{thm:lipschitz-f}.
Since $f^{k+1}=\Psi(f^k)$ for all $k\in\mathbb N$, passing
to the limit gives
$
\Psi(f^\infty)
=
\lim_{k\to\infty}\Psi(f^k)
=
\lim_{k\to\infty}f^{k+1}
=
f^\infty.
$
Hence $f^\infty$ is a Lipschitz continuous BNE.

Finally, we show that the limit $f^\infty$ must be the greatest equilibrium. To see this,
let $\tilde{f} \neq f^\infty$ be a BNE. Since $f^0 \succeq \tilde{f}$,
by applying the operator $\Psi(\cdot)$ on both sides 
of the inequality iteratively over all $k$ and using the fact that
$\Psi(\tilde{f}) = \tilde{f}$, we 
arrive at $f^k \succeq \tilde{f}$ for all $k \in \mathbb{N}$. Letting $k\to\infty$, we 
obtain $f^\infty \succeq \tilde{f}$. Therefore, $f^\infty$ is the greatest equilibrium.

\underline{Part (ii)}.
Proposition \ref{prop: monotonicity-f}(ii) guarantees that, for each $f_{-i} \in \F_{-i}^+$, the optimal response function $\theta_i \mapsto \A_i^*(f_{-i},\theta_i)$ is Lipschitz continuous and increasing on $\Theta_i$, and hence $\A_i^*(f_{-i},\cdot)\in \C_i^+$.
Therefore, the operator $\Psi : \C^+ \to \C^+$.
Since $f^0 = \left(
\vee \A_1,\dots, \vee \A_n \right)
\in \C^+$, it implies that
$f^1 = \Psi(f^0) \in \C^+$. 
By the same induction scheme as in Part (i), we can derive the greatest Lipschitz continuous equilibrium $f^\infty$.
For all $\theta \geq \theta' \in \Theta$ and 
$k \in \mathbb{N}$,
we have 
$f^k(\theta) \geq f^k(\theta')$,
then it follows that
$
f^\infty(\theta) =
\lim_{k \to \infty} f^k(\theta) \geq 
\lim_{k \to \infty}
f^k(\theta')
= f^\infty(\theta') 
$.
Therefore,
the greatest equilibrium
$f^\infty$ is increasing over $\Theta$.

\underline{Part (iii)}.
When $\rho(M) < 1$,  the uniqueness result follows from Theorem~\ref{thm:exist-unique-cbne}(ii) directly.
$\hfill \square$

In Theorem~\ref{thm:exist-unique-CMBNE}, the order-theoretic monotone iteration used in Parts~(i)--(ii) follows the classical construction in monotone Bayesian games. The contribution here is not the monotone iteration itself, but its combination with the Lipschitz continuity and our contraction-based uniqueness results established in Theorem~\ref{thm:BNE-contraction} and Proposition~\ref{prop:rhoM<1}.

\section{Stability of CBNE}
\label{se:stability}

In this section, we investigate stability of continuous Bayesian Nash equilibrium against the perturbation of $\eta$, the joint distribution of $\theta$. 
This is particularly relevant when the equilibrium being studied is based on the present data, and we would like to know whether such an equilibrium is stable when the data is shifted for various reasons.
For example, in the price competition model, the equilibrium strategy $f^*$ describes each firm's pricing rule where $\theta_i$ represents firm $i$'s marginal cost.
Changes in market conditions, geopolitical risks, or policy environments may alter firms' beliefs about their rivals' private costs, even before such changes are fully reflected in historical data.
Such a shift in the information environment can, in turn, be viewed as a perturbation of the underlying type distribution, which may lead to changes in pricing equilibrium, e.g., moving upward or downward with respect to belief shifts and the magnitude of the shift.
The former is related to the comparative-statics literature on Bayesian games with strategic complementarities (see e.g.~\cite{vives2007monotone}). Our focus in this section will be on the latter, the magnitude of the shift.

Let $\mathscr{P}(\Theta)$ denote the set of all probability distributions over $\Theta$ and 
$\eta \in \mathscr{P}(\Theta)$ be the true probability distribution of $\theta$
as discussed in the preceding sections. 
We consider a perturbation $\mu \in \mathscr{P}(\Theta)$ 
and the conditional distribution $\mu_i(\cdot\mid\theta_i)$ which 
represents player $i$'s perturbed posterior belief about rivals' types after observing its own type $\theta_i$.
In this section, we restrict the perturbed 
distributions to the set 
\begin{equation}
\label{def:M-Theta}
\M(\Theta) :=
\bigg\{
\mu \in \mathscr P(\Theta):
\text{for each } S\in\mathscr B(\Theta_{-i}),\ 
\mu_i(S\mid\theta_i)\text{ is measurable in } \theta_i,  \forall\, i\in N
\bigg\},
\end{equation}
which will ensure that the objective function of each player under the perturbed distribution is measurable. 
Compared to 
the true distribution $\eta$, the requirement 
for the perturbed distribution $\mu$ is 
weaker, i.e., 
Lipschitz continuity
of the conditional distributions with respect to the player's own type is not needed.
For each $\mu \in {\cal M}(\Theta)$ and 
$i\in N$,
define the expected utility functions
\begin{subequations}
\bgeqn 
\label{eq:perturbation-1}
\phi_{i, \eta}(a_i, f_{-i}, \theta_i) := \int_{\Theta_{-i}} u_i(a_i, f_{-i}(\theta_{-i}), \theta) d \eta_i \left(\theta_{-i} | \theta_i\right), \
\forall \theta_i \in \Theta_i,\\
\label{eq:perturbation-2}
\phi_{i, \mu}(a_i, f_{-i}, \theta_i) := \int_{\Theta_{-i}} u_i(a_i, f_{-i}(\theta_{-i}), \theta) d \mu_i \left(\theta_{-i} | \theta_i\right), \
\forall \theta_i \in \Theta_i.
 \edeqn  
\end{subequations}
and optimal response functions
\begin{equation}
\label{eq:perturbation-3}
\A_{i,\eta}^*(f_{-i}, \theta_i):= \mathop{\arg \max}_{a_i \in \A_i} 
\ \phi_{i, \eta}(a_i, f_{-i}, \theta_i),\
\A_{i,\mu}^*(f_{-i}, \theta_i):= \mathop{\arg \max}_{a_i \in \A_i} 
\ \phi_{i, \mu}(a_i, f_{-i}, \theta_i),
\end{equation}
in conjunction with their true counterparts.
A prerequisite for the stability analysis is that for every 
$\mu \in \M(\Theta)$,
\begin{equation}
\label{eq:gradient_phi_ai_mu}
\nabla_{a_i}\phi_{i,\mu}(a_i,f_{-i},\theta_i)
=
\int_{\Theta_{-i}}
\nabla_{a_i}u_i(a_i,f_{-i}(\theta_{-i}),\theta_i,\theta_{-i})
\,d\mu_i(\theta_{-i}\mid\theta_i).
\end{equation}
The equality is 
down to the interchangeability between differentiation and the integration, which is guaranteed  
by Assumption~\ref{assu:u_i}(a) and
Assumption~\ref{ass:lipschitz}(a).

Assuming that $\A_{i,\mu}^*(f_{-i}, \theta_i)$ and $\A_{i,\eta}^*(f_{-i}, \theta_i)$ are singleton, we quantify 
the difference between the corresponding 
CBNEs in terms of the Kantorovich distance between the probability distributions
$\eta_i(\cdot\mid\theta_i)$ and $ \mu_i(\cdot\mid\theta_i) $ for $i\in N$.
We need the following technical assumption.

\begin{assumption}
\label{ass:lipschitz-3}
Consider the BNE model (\ref{eq:BNE}). 
For $i \in N$, 
$\nabla_{a_i} u_i(a_i, a_{-i}, \theta_i, \theta_{-i})$ is Lipschitz continuous in $\theta_{-i}$ uniformly with respect to
$a$ and $\theta_i$ with constant modulus $\varrho_i$, i.e., 
\begin{equation}
\label{eq:lipschitz-3}
\|
\nabla_{a_i}u_i(a_i, a_{-i}, \theta_i, \theta_{-i}') - 
\nabla_{a_i}u_i(a_i, a_{-i}, \theta_i, \theta_{-i}'')
\|   \leq
\varrho_i\|\theta_{-i}' - \theta_{-i}'' \|, \;\; \forall \theta_{-i}', \theta_{-i}'' \in \Theta_{-i}.
\end{equation}
\end{assumption}

The condition requires each player's marginal utility function to be Lipschitz continuous w.r.t.~change of rivals' type parameters. The condition is not demanding since type parameters of rivals may even 
not directly affect a player's utility function or its marginal in some practical application (in which case $\varrho_i=0$).

\begin{theorem}[Pointwise stability of CBNE]
\label{thm:belief-perturb}
Suppose: 
(a) Assumptions~\ref{assu:u_i} - \ref{ass:lipschitz} and \ref{ass:lipschitz-3} 
hold;
(b) there exists a vector $w\in \R_{++}^n$ such that $\alpha=\|D^{-1}MD\|_\infty<1$, 
where $D=\mathrm{diag}(w)$.
Let $f^*_\eta$ be the true CBNE and $f^*_\mu$ be the unique BNE under perturbation $\mu\in {\cal M}(\Theta)$. Then
\begin{equation}
\label{eq:f-stab-infty}
\| f^*_\eta -  f^*_\mu\|_{\infty,w}
:=
\max_{i \in N}
\left(
\frac{1}{w_i}\sup_{\theta_i \in \Theta_i}\| f^*_{i,\eta}(\theta_i) - f^*_{i, \mu}(\theta_i) \|
\right)
\leq
\max_{i \in N}
\sup_{\theta_i \in \Theta_i}
\frac{\tau_i + \varrho_i}{w_i\sigma_i (1-\alpha)}
\dd_{\mathrm{K}} \big(\eta_i(\cdot\mid\theta_i), \mu_i(\cdot\mid\theta_i) \big),
\end{equation}
where 
the Kantorovich (Wasserstein–1) distance between $\mu,\nu \in \mathscr P(\Xi)$ is defined by
$$
\dd_{\mathrm{K}}(\mu,\nu) :=
\sup\left\{
\left|\int_\Xi h\,d\mu-\int_\Xi h\,d\nu\right|
:\ h:\Xi\to \R,\ \text{Lip}(h) \leq 1
\right\},
$$ 
and 
$\tau_i:=
\left(\sum_{j \neq i} \tau_{ij}^2 \left(\frac{\kappa_j}{\sigma_j}\right)^2 \right)^{{1/2}}$.
\end{theorem}

Before providing the proof, it might be helpful to make some comments about the conditions of the theorem. 
Under conditions (a) and (b), existence and uniqueness of the perturbed equilibrium $f_\mu^*$ are guaranteed by the same contraction argument as in Theorem~\ref{thm:BNE-contraction}.
Moreover, for each 
observed type $\theta_i$, the conditional distributions $\eta_i(\cdot\mid\theta_i)$ and $\mu_i(\cdot\mid\theta_i)$ represent player $i$'s beliefs about rivals' types before and after the perturbation, the Kantorovich distance measures how much these two conditional beliefs differ in a geometric sense, namely, how much probability mass must be reallocated across the rivals' type space in order to move from one belief to the other.

\noindent
\textbf{Proof.}
We use Lemma~\ref{Lem:Lip-slv-mon} to prove the result.
Under Assumption~\ref{ass:lipschitz-3} and the conditions in Theorem~\ref{thm:exist-unique-cbne}, 
$
-\nabla_{a_i} \phi_i\left(a_i, f_{-i}, \theta_i\right)$ 
is strongly monotone in $a_i$ uniformly with respect to 
$f_{-i}, \theta_i$ 
and conditional probability distribution
$\eta_i(\cdot \mid \theta_i)$. 
We consider $f \in \Psi(\F)$ (w.r.t.~$\eta$)
by fixing $f_{-i}, \theta_i$ and treating 
$\eta_i(\cdot \mid \theta_i)$ as a varying parameter, we obtain by virtue of Lemma~\ref{Lem:Lip-slv-mon} that for any $\mu \in \M(\Theta)$,
\begin{align}
\label{eq:eta-mu}
\nonumber
& \|\A_{i,\eta}^*(f_{-i}, \theta_i) - \A_{i,\mu}^*(f_{-i}, \theta_i) \| \\
\nonumber
& \leq  \frac{1}{\sigma_i}
\|
\nabla_{a_i} \phi_{i,\eta}(\A_{i,\eta}^*(f_{-i}, \theta_i), f_{-i}, \theta_i)
- \nabla_{a_i} \phi_{i,\mu}(\A_{i,\eta}^*(f_{-i}, \theta_i), f_{-i}, \theta_i)
\| \\ 
& = 
\frac{1}{\sigma_i}
\left\|
\int_{\Theta_{-i}}
\nabla_{a_i}
u_i \big(\A_{i,\eta}^*(f_{-i}, \theta_i), f_{-i}(\theta_{-i}), \theta \big)
\big(
d\eta_i(\theta_{-i} | \theta_i) - d\mu_i(\theta_{-i} | \theta_i)
\big) \right\|. 
\end{align}
Observe that $\forall \theta_{-i}', \theta_{-i}'' \in \Theta_{-i}$,
\begin{align*}
& \left \|
\nabla_{a_i}
u_i \big(\A_{i,\eta}^*(f_{-i}, \theta_i), f_{-i}(\theta_{-i}'), \theta_i, \theta_{-i}' \big)
-
\nabla_{a_i}
u_i \big(\A_{i,\eta}^*(f_{-i}, \theta_i), f_{-i}(\theta_{-i}''), \theta_i, \theta_{-i}'' \big) 
\right \| \\
& \leq
\left \|
\nabla_{a_i}
u_i \big(\A_{i,\eta}^*(f_{-i}, \theta_i), f_{-i}(\theta_{-i}'), \theta_i, \theta_{-i}' \big)
-
\nabla_{a_i}
u_i \big(\A_{i,\eta}^*(f_{-i}, \theta_i), f_{-i}(\theta_{-i}''), \theta_i, \theta_{-i}' \big) 
\right \| \\
& \quad + 
\left \|
\nabla_{a_i}
u_i \big(\A_{i,\eta}^*(f_{-i}, \theta_i), f_{-i}(\theta_{-i}''), \theta_i, \theta_{-i}' \big)
-
\nabla_{a_i}
u_i \big(\A_{i,\eta}^*(f_{-i}, \theta_i), f_{-i}(\theta_{-i}''), \theta_i, \theta_{-i}'' \big) 
\right \| \\
& \leq
\sum_{j \neq i} \tau_{ij}
\|f_j(\theta_j') - f_j(\theta_j'') \| 
+ \varrho_i \|\theta'_{-i} - \theta_{-i}'' \| 
\quad \text{(by \eqref{eq:lipschitz-2} and \eqref{eq:lipschitz-3})}
\\
& \leq 
\sum_{j \neq i} \tau_{ij}\frac{\kappa_j}{\sigma_j}
\|\theta_j' - \theta_j'' \| 
+ \varrho_i \|\theta'_{-i} - \theta_{-i}'' \| 
\quad \text{(by \eqref{eq:lipschitz-theta})}
\\
& \leq 
\left(\sum_{j \neq i} 
\tau_{ij}^2  \left(\frac{\kappa_j}{\sigma_j}\right)^2 \right)^{1/2}
\|\theta_{-i}' - \theta_{-i}'' \|
+ \varrho_i \|\theta'_{-i} - \theta_{-i}'' \| 
\quad \text{(by Cauchy--Schwarz inequality)}
\\
& =
(\tau_i + \varrho_i) 
\|\theta'_{-i} - \theta_{-i}'' \|,
\end{align*}
where the third inequality follows from the fact that
$f \in \Psi(\F)$ (w.r.t. $\eta$).
Indeed, there exists $f' \in\F$ such that
$f= \left( \A_{i,\eta}^*(f_{-i}', \cdot) \right)_{i \in N}.$
This shows that  
$\nabla_{a_i}
u_i \big(\A_{i,\eta}^*(f_{-i}, \theta_i), f_{-i}(\theta_{-i}), \theta_i, \theta_{-i} \big)$
is Lipschitz continuous in $\theta_{-i}$ over $\Theta_{-i}$ with modulus $(\tau_i + \varrho_i)$.
By the definition of Kantorovich metric
\bgeq 
\left\|
\int_{\Theta_{-i}}
\nabla_{a_i}
u_i \big(\A_{i,\eta}^*(f_{-i}, \theta_i), f_{-i}(\theta_{-i}), \theta \big)
\big(
d\eta_i(\theta_{-i} | \theta_i) - d\mu_i(\theta_{-i} | \theta_i)
\big) \right\| \leq 
(\tau_i + \varrho_i)
\dd_{\mathrm{K}} \big(\eta_i(\cdot\mid\theta_i), \mu_i(\cdot\mid\theta_i) \big).
\edeq
Combining this inequality with \eqref{eq:eta-mu}, we obtain
\begin{equation}
\label{eq:optimal-response-stability}
\|\A_{i,\eta}^*(f_{-i}, \theta_i) - \A_{i,\mu}^*(f_{-i}, \theta_i) \|
\leq
\frac{\tau_i + \varrho_i}{\sigma_i }
\dd_{\mathrm{K}} \big(\eta_i(\cdot\mid\theta_i), \mu_i(\cdot\mid\theta_i) \big).
\end{equation}
We are now ready to quantify the difference between $f^*_\eta$ and $  f^*_\mu$. We need the mapping $\Psi$ defined in \eqref{eq:PSI}.
Note that since the function $\phi_i$ defined in \eqref{eq:def-phi_i-BNE} depends on the underlying distribution $\eta$, then $\Psi$ also depends on the distribution. 
Thus we write $\Psi_\eta$ and $\Psi_\mu$ respectively when the joint distribution of $\theta$ is perturbed from $\eta$ to $\mu$. By the triangle inequality and \eqref{eq:optimal-response-stability}, we obtain
\begin{align}
\label{eq:f-stb}
\nonumber
\| f^*_\eta -  f^*_\mu\|_{\infty,w}
& =
\| \Psi_\eta(f^*_\eta) - \Psi_\mu(f^*_\mu)\|_{\infty,w} 
\leq
\| \Psi_\eta(f^*_\eta) - \Psi_\mu(f^*_\eta)\|_{\infty,w}
+
\| \Psi_\mu(f^*_\eta) - \Psi_\mu(f^*_\mu)\|_{\infty,w} 
\\ \nonumber
& \leq  
\max_{i \in N}
\frac{1}{w_i}
\|\A_{i,\eta}^*(f_{-i,\eta}^*, \cdot) - \A_{i,\mu}^*(f_{-i,\eta}^*, \cdot) \|_{\infty}
+
\alpha \|f^*_\eta - f^*_\mu\|_{\infty,w}
\\ 
& \leq
\alpha \|f^*_\eta - f^*_\mu\|_{\infty,w}
+ 
\max_{i \in N}
\sup_{\theta_i \in \Theta_i}
\frac{\tau_i + \varrho_i}{w_i \sigma_i }
\dd_{\mathrm{K}} \big(\eta_i(\cdot\mid\theta_i), \mu_i(\cdot\mid\theta_i) \big),
\end{align}
where $\alpha=\|D^{-1}MD\|_\infty$ is defined in 
\eqref{eq:alpha}. 
A reorganization of the inequality above gives rise to
\eqref{eq:f-stab-infty}.
$\hfill$ $\square$

Inequality \eqref{eq:f-stab-infty} provides a pointwise (scenario-based) bound on the difference between the BNEs before and after perturbation of the probability distribution of $\theta$.
The error bound is controlled linearly by the Kantorovich distance between the two conditional distributions, which means 
that the difference between 
the optimal response functions 
of the two BNEs can be small 
in each scenario of $\theta$ 
when the perturbed distribution is sufficiently close to the true distribution.
The deficiency of this result is that the error bound is determined by the largest deviation of 
one of the players' conditional distributions. 
This raises the question as to whether we can quantify the difference by ${\cal L}^{p}$-norm. To address the question, we need an additional condition on the marginal distributions of the type parameters.

\begin{assumption}[Equal own-type marginals]
\label{ass:marginal-density}
Let $\mu\in \M(\Theta)$. 
For each player $i \in N$, the
marginal distribution of the own type $\theta_i$ under $\mu$ coincides with that under $\eta$, that is,
\begin{equation}
\label{eq:marginal-density=}
\eta_i(S) = \mu_i(S),\quad \text{for every Borel set } S \in \mathscr{B}(\Theta_i).
\end{equation}
\end{assumption}
Assumption~\ref{ass:marginal-density} 
ensures that the joint distributions $\eta$ and $\mu$ share the same marginal, which means that each player's belief about its own type is not affected by the perturbation, while its beliefs about the distributions of rivals' types and the dependence structure between  $\theta_i$ and $\theta_{-i}$ may change.
This reflects the fact that, in practice, a player typically has more information about its own type parameter than about those of its rivals.
Under the assumption, we can use 
$\mathcal{L}^p(\eta_i)$-norm to assess
the optimal response gap $\Delta_i(f_{-i}, \theta_i):=\A^*_{i,\eta}(f_{-i},\theta_i)-\A^*_{i,\mu}(f_{-i},\theta_i)$, for each fixed $f_{-i} \in \F_{-i}$ with
$$
\| \Delta_i(f_{-i}) \|_{\mathcal{L}^p(\eta_i)}
:=
\begin{cases}
\left(\displaystyle\int_{\Theta_i}\|\Delta_i(f_{-i}, \theta_i)\|^p d\eta_i(\theta_i) \right)^{1/p},
& \text{for} \; 1\leq p<\infty,\\[0.6em]
\underset{\theta_i\in\Theta_i}{\mathrm{ess\,sup}}\ \|\Delta_i(f_{-i}, \theta_i)\|,
& \text{for}\; p=\infty,
\end{cases}
$$
which averages type–wise deviations according to its own type distribution.
By the pointwise sensitivity bound in \eqref{eq:optimal-response-stability},
$
\|\Delta_i(f_{-i}, \theta_i)\|
\ \leq\
\frac{\tau_i+\varrho_i}{\sigma_i}
\dd_{\mathrm{K}}\left(\eta_i(\cdot\mid\theta_i),\mu_i(\cdot\mid\theta_i)\right).
$
For fixed $1\leq p < \infty$ and $f_{-i} \in \F_{-i}$, we may take the $\mathcal{L}^p(\eta_i)$–norm to obtain 
\begin{equation}
\label{eq:lp-conditionalW1}
\|\Delta_i(f_{-i})\|_{\mathcal{L}^p(\eta_i)}
\leq
\frac{\tau_i+\varrho_i}{\sigma_i}
\left(\int_{\Theta_i}
\left(
\dd_{\mathrm{K}}\big(\eta_i(\cdot\mid\theta_i),\mu_i(\cdot\mid\theta_i)\big) 
\right)^p
d\eta_i(\theta_i)\right)^{1/p}.
\end{equation}
When $p=1$, the term on the right-hand side of \eqref{eq:lp-conditionalW1} is an averaged conditional Kantorovich distance, which is closely related to the conditional Wasserstein distance studied by Chemseddine et al.~\cite{chemseddine2025conditional}. 
They establish a bound of the form
\begin{equation}
\label{eq:conditional-wasser}
\dd_{\mathrm{K}}(\eta, \mu) \leq 
\left(\int_{\Theta_i}
\dd_{\mathrm{K}}\big(\eta_i(\cdot\mid\theta_i),\mu_i(\cdot\mid\theta_i)\big) d\eta_i(\theta_i)\right).
\end{equation}
Unfortunately, we cannot use 
inequality \eqref{eq:conditional-wasser}
as we want to bound
$\dd_{\mathrm{K}}\big(\eta_i(\cdot\mid\theta_i),\mu_i(\cdot\mid\theta_i)\big) $
in terms of $\dd_{\mathrm{K}}(\eta, \mu)$ in \eqref{eq:f-stab-infty}.  
This motivates us to adopt the relationship 
between the Kantorovich metric and the KL-divergence 
bridged by the total variation metric when $\Theta$ is bounded, as established in Gibbs and Su~\cite{gibbs2002choosing}.

\begin{theorem}[Stability and sensitivity of CBNE in $\eta$]
\label{thm:belief-perturb-2}
Assume the setting and conditions of Theorem~\ref{thm:belief-perturb}.
For $i \in N$, let  $\inmat{Diam}(\Theta_{-i})$ be the diameter of set $\Theta_{-i}$. 
Then the following assertions hold.
\begin{itemize}
\item[(i)] For any perturbation distribution $\mu\in {\cal M}(\Theta)$ 
satisfying $\eta \ll \mu$ and Assumption~\ref{ass:marginal-density}, 
\begin{equation}
\label{eq:f-stb-L_2}
\|f^*_\eta - f^*_\mu \|_{\infty, w, \mathcal{L}^2(\eta)}
\leq
\max_{i \in N}
\left(\frac{\tau_i + \varrho_i}{w_i \sigma_i(1-\alpha) } \inmat{Diam}(\Theta_{-i})\right)
\sqrt{\frac{\dd_{\mathrm{KL}}(\eta\ \| \ \mu)}{2}},
\end{equation}
where 
$$
\dd_{\mathrm{KL}}(\mu \ \| \ \nu)
:=
\begin{cases}
\displaystyle
\int_\Xi \log\left(\frac{d\mu}{d\nu}\right)\,d\mu,
& \text{if } \mu \ll \nu,\\[1ex]
+\infty,
& \text{otherwise,}
\end{cases}
$$
denotes 
the Kullback–Leibler (KL) divergence 
of $\mu$ from $\nu$ for $\mu,\nu \in \mathscr P(\Xi)$.

\item[(ii)] 
Let $\hat{\eta} \in \M(\Theta)$
satisfying $\eta \ll \hat{\eta}$ and consider perturbation distribution $\mu :=  (1 - \epsilon) \eta + \epsilon \hat{\eta}$. 
If $\hat{\eta}$ satisfies Assumption~\ref{ass:marginal-density}, then
\begin{equation}
\label{eq:f-stv-L_2}
\limsup_{\epsilon\downarrow0} \frac{\|f^*_\eta - f^*_\mu \|_{\infty, w, \mathcal{L}^2(\eta)}}{\epsilon}
\leq
\max_{i \in N}
\left(\frac{\tau_i + \varrho_i}{w_i \sigma_i(1-\alpha) } \inmat{Diam}(\Theta_{-i})\right)
\sqrt{\frac{\dd_{\mathrm{KL}}(\eta \ \| \  \hat{\eta})}{2}}.
\end{equation}

\end{itemize}
\end{theorem}

It might be helpful 
 to comment on the condition $\eta \ll \mu$ and the associated error bounds.
First, by requiring $\eta$ to be absolutely continuous with respect to $\mu$, we interpret
$\mu$ in practice as the target distribution under which the BNE model will be applied in
the future, while $\eta$ represents the current reference distribution under which the BNE
model is analyzed. The error bound in~\eqref{eq:f-stab-infty} is expressed in terms of
conditional distributions, which is well suited for interim equilibrium analysis, where
each player has complete information about its own type. In contrast, 
the error bound
in~\eqref{eq:f-stb-L_2} controls the deviation of the joint distribution $\eta$ from $\mu$,
making it more suitable for ex ante equilibrium analysis, where each player is uncertain
about its own type.
In part (ii), by setting $\mu :=  (1 - \epsilon) \eta + \epsilon \hat{\eta}$, we mean that the random data of the type parameters
are 
obtained from two different sources, one generated by
$\eta$ with probability $(1-\epsilon)$ and 
the other generated by 
$\hat{\eta}$ with probability $\epsilon$.

\noindent
\textbf{Proof.}
\underline{Part (i)}. 
For any $f \in \Psi_\eta(\F)$ and fixed $\theta_i\in\Theta_i$, we have by inequality~\eqref{eq:optimal-response-stability}
\begin{align}
\label{eq:gap-square}
\nonumber
& \|\A_{i,\eta}^*(f_{-i}, \theta_i) - \A_{i,\mu}^*(f_{-i}, \theta_i) \|^2 \\
& \leq \nonumber
\left(\frac{ \tau_i + \varrho_i}{\sigma_i }\right)^2
\left(\dd_{\mathrm{K}} \big(\eta_i(\cdot\mid\theta_i), \mu_i(\cdot\mid\theta_i) \big) \right)^2 
\leq \nonumber
\left(\frac{ \tau_i + \varrho_i}{\sigma_i }\right)^2
\left(\inmat{Diam}(\Theta_{-i}) \dd_{\mathrm{TV}} \big(\eta_i(\cdot\mid\theta_i), \mu_i(\cdot\mid\theta_i) \big) \right)^2 \\
& \leq
\left(\frac{\tau_i + \varrho_i}{\sigma_i } \inmat{Diam}(\Theta_{-i})\right)^2 
\frac{\dd_{\mathrm{KL}} \big(\eta_i(\cdot\mid\theta_i) \ \| \  \mu_i(\cdot\mid\theta_i) \big)}{2},
\end{align}
where $\dd_{\mathrm{TV}}$ denotes the total variation distance, the second inequality is based on 
Lemma~\ref{metric} in the Appendix, 
and the last inequality is from Pinsker's inequality directly.
By the definition of the KL-divergence and tower property, 
\begin{align}
\label{eq:expect-KL}
\nonumber
\int_{\Theta_i}
\dd_{\mathrm{KL}} \big(\eta_i(\cdot\mid\theta_i) \ \| \  \mu_i(\cdot\mid\theta_i) \big)
d \eta_i(\theta_i)
& 
=
\int_{\Theta_i}
\left(
\int_{\Theta_{-i}}
\log\left(
\frac{d\eta_i(\cdot\mid\theta_i)}{d\mu_i(\cdot\mid\theta_i)}(\theta_{-i})
\right)
d\eta_i(\theta_{-i}\mid\theta_i)
\right)
d\eta_i(\theta_i) 
\\
\nonumber
& 
=
\int_{\Theta}
\log\left(
\frac{d\eta}{d\mu}(\theta)  
\bigg /
\frac{d\eta_i}{d\mu_i}(\theta_i)
\right)
d\eta(\theta) 
\\
\nonumber
& 
=
\int_{\Theta}
\log\left(
\frac{d\eta}{d\mu}(\theta)
\right)
d\eta(\theta)
-
\int_{\Theta_i}
\log\left(
\frac{d\eta_i}{d\mu_i}(\theta_i)
\right)
d\eta_i(\theta_i) 
\\ & 
=
\dd_{\mathrm{KL}}(\eta \ \| \ \mu)
-
\dd_{\mathrm{KL}}(\eta_i \ \| \ \mu_i),
\end{align}
which, together with the nonnegativity of the KL-divergence, implies that
\begin{equation}
\label{eq:expect-KL-ineq}
\int_{\Theta_i}
\dd_{\mathrm{KL}} \big(\eta_i(\cdot\mid\theta_i) \ \| \ \mu_i(\cdot\mid\theta_i) \big)
\, d \eta_i(\theta_i)
\leq
\dd_{\mathrm{KL}}(\eta \ \| \ \mu),
\end{equation}
where the equality holds if and only if 
$\dd_{\mathrm{KL}}(\eta_i \ \| \ \mu_i) = 0$.
The latter is implied by 
$\eta_i = \mu_i$.
Combining \eqref{eq:gap-square}-\eqref{eq:expect-KL-ineq}, we obtain
\begin{align}
\label{eq:L2-norm-KL}
\nonumber
\|\A_{i,\eta}^*(f_{-i}, \cdot) - \A_{i,\mu}^*(f_{-i}, \cdot) \|_{\mathcal{L}^2(\eta_i)}^2
& =
\int_{\Theta_i} \|\A_{i,\eta}^*(f_{-i}, \theta_i) - \A_{i,\mu}^*(f_{-i}, \theta_i) \|^2 d\eta_i(\theta_i)\\
& \leq
\left(\frac{\tau_i + \varrho_i}{\sigma_i } \inmat{Diam}(\Theta_{-i})\right)^2  
\frac{\dd_{\mathrm{KL}}(\eta \ \| \  \mu)}{2}.
\end{align}
Since $f_\eta^*=\Psi_\eta(f_\eta^*) \in \Psi_\eta(\F)$ and $f_\mu^*=\Psi_\mu(f_\mu^*) \in \Psi_\mu(\F)$,
\begin{align}
\nonumber
\|f^*_\eta - f^*_\mu \|_{\infty, w, \mathcal{L}^2(\eta)}
& =
\| \Psi_\eta(f^*_\eta) - \Psi_\mu(f^*_\mu)\|_{\infty, w, \mathcal{L}^2(\eta)}
\\ 
& \leq
\| \Psi_\eta(f^*_\eta) - \Psi_\mu(f^*_\eta)\|_{\infty, w, \mathcal{L}^2(\eta)}
+
\| \Psi_\mu(f^*_\eta) - \Psi_\mu(f^*_\mu)\|_{\infty, w, \mathcal{L}^2(\eta)}. 
\label{eq:thm7-proof-2}
\end{align}
Under Assumption~\ref{ass:marginal-density}, and by the contraction property of $\Psi_\mu$, it follows that
\bgeqn 
\| \Psi_\mu(f^*_\eta) - \Psi_\mu(f^*_\mu)\|_{\infty, w, \mathcal{L}^2(\eta)}
&=&
\| \Psi_\mu(f^*_\eta) - \Psi_\mu(f^*_\mu)\|_{\infty, w, \mathcal{L}^2(\mu)}
\leq
\alpha \|f^*_\eta - f^*_\mu\|_{\infty, w, \mathcal{L}^2(\mu)}\nonumber\\
&=& 
\alpha \|f^*_\eta - f^*_\mu\|_{\infty, w, \mathcal{L}^2(\eta)}.
\label{eq:thm7-proof-3}
\edeqn 
Combining \eqref{eq:thm7-proof-2}-\eqref{eq:thm7-proof-3} and using the definition of $\|\cdot\|_{\infty, w, \mathcal{L}^2(\eta)}$-norm, we have
\begin{align}
\label{eq:errorbound}
\|f^*_\eta - f^*_\mu \|_{\infty, w, \mathcal{L}^2(\eta)} \leq
\max_{i \in N}
\frac{1}{w_i}
\|\A_{i,\eta}^*(f_{-i,\eta}^*, \cdot)
-
\A_{i,\mu}^*(f_{-i,\eta}^*, \cdot) \|_{\mathcal{L}^2(\eta_i)}
+
\alpha \|f^*_\eta - f^*_\mu\|_{\infty, w, \mathcal{L}^2(\eta)}.
\end{align}
A reorganization of the inequality 
together with the inequality~\eqref{eq:L2-norm-KL} gives rise to \eqref{eq:f-stb-L_2}.

\underline{Part (ii)}.
By the definition of the Kantorovich metric,
\begin{align}
\label{eq:dK-epsilon}
\dd_{\mathrm{K}} \big(\eta_i(\cdot\mid\theta_i), \mu_i(\cdot\mid\theta_i)\big)
& = 
\sup_{\text{Lip}(g) \leq 1}
\left \vert
\int_{\Theta_{-i}}
g(\theta_{-i})
d\eta_i(\theta_{-i} | \theta_i)
-
\int_{\Theta_{-i}}
g(\theta_{-i})
d\mu_i(\theta_{-i} | \theta_i)
\right \vert  
\end{align}
Under Assumption~\ref{ass:marginal-density},
the own-type marginal distributions coincide, i.e.,
$\mu_i = \eta_i=\hat{\eta}_i$.
Moreover, for any Borel sets $A \in \mathscr{B}(\Theta_i)$ and $B\in \mathscr{B}(\Theta_{-i})$, define
$\tilde{\mu}_i(B\mid\theta_i)
:=
(1-\epsilon)\eta_i(B\mid\theta_i)
+\epsilon\hat{\eta}_i(B\mid\theta_i)$.
Then it follows that
\begin{align*}
\int_A \tilde{\mu}_i(B\mid\theta_i)\,d\mu_i(\theta_i)
&=
(1-\epsilon)\int_A \eta_i(B\mid\theta_i)\,d\eta_i(\theta_i)
+\epsilon\int_A \hat{\eta}_i(B\mid\theta_i)\,d\hat{\eta}_i(\theta_i) \\
&=
(1-\epsilon)\eta(A\times B)+\epsilon\hat{\eta}(A\times B) 
= \mu(A\times B),
\end{align*}
which implies that $\tilde{\mu}_i(\cdot\mid\theta_i)$ is a regular conditional distribution of $\mu$
given $\theta_i$. 
Therefore, the equality
$\mu_i(\cdot\mid\theta_i)
= \tilde{\mu}_i (\cdot\mid\theta_i)$ holds $\mu_i\text{-a.e. }\theta_i\in\Theta_i$.
Since $\mu_i=\eta_i$, for $\eta_i$-a.e. $\theta_i\in\Theta_i$, we have
\begin{align}
\label{eq:dK-epsilon-1}
\nonumber
\dd_{\mathrm{K}} \big(\eta_i(\cdot\mid\theta_i), \mu_i(\cdot\mid\theta_i)\big)
& =
\sup_{\text{Lip}(g) \leq 1}
\left|
\epsilon
\int_{\Theta_{-i}} g(\theta_{-i})\,
d\big(\eta_i(\cdot\mid\theta_i)-\hat{\eta}_i(\cdot\mid\theta_i)\big)(\theta_{-i})
\right| \\
& =
\epsilon\,
\dd_{\mathrm{K}} \big(\eta_i(\cdot\mid\theta_i), \hat{\eta}_i(\cdot\mid\theta_i)\big).
\end{align}
Using the same argument as in \eqref{eq:gap-square} together with \eqref{eq:dK-epsilon-1}, we obtain for $\eta_i$-a.e. $\theta_i\in\Theta_i$,
\begin{align}
\label{eq:final}
\nonumber
& \|\A_{i,\eta}^*(f_{-i}, \theta_i) - \A_{i,\mu}^*(f_{-i}, \theta_i) \|^2 
\\ \nonumber
& \leq 
\left(\frac{\tau_i + \varrho_i}{\sigma_i }\right)^2
\left(\epsilon \dd_{\mathrm{K}} \big(\eta_i(\cdot\mid\theta_i), \hat{\eta}_i(\cdot \mid \theta_i) \big) \right)^2  
\leq
\left(\frac{\tau_i + \varrho_i}{\sigma_i }\right)^2
\left(\epsilon \inmat{Diam}(\Theta_{-i}) \dd_{\mathrm{TV}} \big(\eta_i(\cdot\mid\theta_i), \hat{\eta}_i(\cdot\mid\theta_i) \big) \right)^2 \\
& \leq
\left(\frac{\tau_i + \varrho_i}{\sigma_i } \epsilon \inmat{Diam}(\Theta_{-i})\right)^2 
\frac{\dd_{\mathrm{KL}} \big(\eta_i(\cdot\mid\theta_i) \ \| \  \hat{\eta}_i(\cdot\mid\theta_i) \big)}{2}.
\end{align}
It follows from \eqref{eq:errorbound} that
\begin{align}
\nonumber
\limsup_{\epsilon\downarrow0} \frac{\|f^*_\eta - f^*_\mu \|_{\infty, w, \mathcal{L}^2(\eta)}}{\epsilon}
& \leq
\limsup_{\epsilon\downarrow0}
\frac{\max\limits_{i \in N} \frac{1}{w_i}
\|\A_{i,\eta}^*(f_{-i, \eta}^*, \cdot) - \A_{i,\mu}^*(f_{-i, \eta}^*, \cdot) \|_{\mathcal{L}^2(\eta_i)}}{\epsilon (1-\alpha)}\\
& =  \nonumber
\limsup_{\epsilon\downarrow0}
\frac{\max\limits_{i \in N} \frac{1}{w_i}
\sqrt{\int_{\Theta_i} \|\A_{i,\eta}^*(f_{-i, \eta}^*, \theta_i) - \A_{i,\mu}^*(f_{-i, \eta}^*, \theta_i) \|^2 d\eta_i(\theta_i)}}{\epsilon (1-\alpha)}
 \\
& \leq
\max_{i \in N}
\left(\frac{\tau_i + \varrho_i}{w_i \sigma_i (1 - \alpha)} \inmat{Diam}(\Theta_{-i})\right)
\sqrt{\frac{\dd_{\mathrm{KL}}(\eta \ \| \  \hat{\eta})}{2}},
\end{align}
where the last inequality follows by integrating \eqref{eq:final} with respect to $\eta_i$ and then applying \eqref{eq:expect-KL-ineq}.
$\hfill \square$

To facilitate understanding of the theorem, we 
use a simple example to illustrate.

\begin{example}[Symmetric two-player Bayesian Cournot game]
\label{example:2p_sym_BCG}
Consider a symmetric Bayesian Cournot game with two players.  
For each player $i=1,2$, let
$$
u_i(a_i,a_{-i},\theta) 
= a_i \, p(a_1,a_2) - c_i(a_i,\theta_i),
$$
where $p(a_1,a_2) = \alpha - \beta(a_1+a_2), \ \alpha,\beta>0,$ is the inverse demand function
and $c_i(a_i,\theta_i) = \theta_i a_i + \frac{c}{2}a_i^2,\ c>0$, is the cost function.
The action and type spaces are
$\A_i = [0,\frac{\alpha}{\beta}], \ 
\Theta_i = [0,1], i=1,2,$
and the type profile $\theta=(\theta_1,\theta_2)$ is drawn from a $\rho$-parametric distribution $\eta_\rho$ supported on $\Theta_1\times \Theta_2$ with density 
$q_\rho(\theta_1, \theta_2) = 1 + \rho(2\theta_1 - 1)(2\theta_2 - 1), -1 <\rho < 1$. 
Then 
$$
\phi_i(a_i,f_{-i},\theta_i) 
= (\alpha-\theta_i)a_i - \Big(\beta+\frac{c}{2}\Big)a_i^2 
-  \beta a_i \int_{\Theta_{-i}} 
f_{-i}(\theta_{-i})
q_i(\theta_{-i}|\theta_i)\, d\theta_{-i}.
$$
We can obtain a closed form of CBNE with 
$$
f_{i, \rho}^*(\theta_i) = 
\frac{\alpha}{3\beta + c} +
\frac{\beta(3 - \rho)}{2(3\beta + c)(6\beta + 3c + \beta \rho)} - \frac{\theta_i}{2\beta + c + \frac{\beta \rho}{3}},
\text{ for } i = 1,2.
$$
The equilibrium is stable w.r.t.~variation of parameter $\rho$ in that for any $\rho_1, \rho_2 \in (-1, 1)$, we have
$$
|f_{i, \rho_1}^*(\theta_i) - f_{i, \rho_2}^*(\theta_i)|
\leq
\left|
\frac{
3\beta(\rho_1-\rho_2)(2\theta_i-1)}
{2(\beta\rho_1+6\beta+3c)(\beta\rho_2+6\beta+3c)}
\right|
\leq
\frac{3\beta}{2(5\beta+3c)^2}\,|\rho_1-\rho_2|.
$$
\end{example}

\section{Concluding Remarks}
\label{se:Concluding remarks}

In this paper, we explore a new way to derive existence and uniqueness of 
BNEs 
by virtue of Banach fixed-point theorem,
which not only simplifies 
the proofs in the literature \cite{guo2021existence,ui2016bayesian},
but also provides new sufficient conditions for existence and uniqueness of BNE
and the business insights behind.
Moreover, our tentative theoretical analysis shows that
the contraction mapping-based 
proof provides a potential avenue, alternative to the existing VI-based approach, 
for computing an approximate BNE in both discontinuous and Lipschitz continuous case.
We also make the first attempt to quantify the impact of perturbation of the joint probability distribution of type parameters on the CBNE by deriving an explicit error bound, which fills a gap in the literature of BNE. 
All of the theoretical results are established with 
bounded support set of type parameters. It remains 
an open and interesting question 
(at least from a theoretical point of view)  whether the results can be extended to the unbounded case. In the stability analysis, we argue that information on the true distribution of type parameters may not be complete and consequently an approximate distribution might have to be used particularly in some data-driven environments. An alternative way to tackle the incomplete information issue is to consider a distributionally robust 
model (e.g.,~\cite{liu2018distributionally,LiuBDRNE2025}). This would require an overhaul from modeling to theoretical analysis and computational methods. In the case when information can be acquired in a learning process, we may also incorporate a Bayesian approach into the BNE model. We leave all these for future research.


\bibliographystyle{plain}

\bibliography{references}

\appendix
\section{Appendix}

\subsection{Stability of strongly monotone parametric variational inequality}

The proofs of some technical results in the main body of the paper 
require frequently 
existence and Lipschitz continuity of the solution to a 
parametric variational inequality where the underlying function is strongly monotone and Lipschitz continuous stated in the next lemma. 

\begin{lemma}
\label{Lem:Lip-slv-mon}
Let $\X \subset \R^m$ be a non-empty, compact, and convex set and $\Y$ be a Polish space. 
Let $\Phi: \X \times \Y \to \R^m$ be a continuous mapping which is $\sigma$-strongly monotone in $x$ uniformly w.r.t.~$y$.
Then, for every $y\in\mathcal Y$, the generalized equation
\begin{equation}
0\in \Phi(x,y)+\N_{\X}(x)
\end{equation}
admits a unique solution $x^*(y)\in\X$, where 
\begin{equation}
\N_{\X} (x) :=
\Big\{ w \in \R^m:\ \langle w, x'- x \rangle \leq 0,\ \forall x' \in \X \Big\} 
\end{equation}
is the normal cone to $\X$ at $x$. 
Moreover, \begin{equation}
\label{eq:Lema-Lip-solu}
\|x^*(y_1) - x^*(y_2) \| \leq 
\frac{1}{\sigma}
\|
\Phi(x^*(y_2), y_1) -\Phi(x^*(y_2), y_2)
\|,\ \forall y_1, y_2\in {\cal Y}.
\end{equation}

\end{lemma}

The proof follows the same line of argument as that of
Dontchev and Rockafellar~\cite[Theorem~2F.6]{rockafellar2009implicit}, which treats the special case
\(\Phi(x,y) := y - g(x)\), where both \(x\) and \(y\) belong to
finite-dimensional spaces and the mapping \(g\) is strongly monotone.
We omit the details.

\subsection{Schauder's fixed point theorem and Arzel\'a-Ascoli theorem}

Let $\Theta$ be a nonempty compact subset of $\R^d$ and let
$Z:=C(\Theta;\R^z)$ be the Banach space of continuous functions from $\Theta$ to $\R^z$ equipped with the infinity norm
$\|f\|_\infty:=\sup_{\theta \in \Theta}\| f(\theta) \|$.
Let $T: Z \rightarrow Z$ be an operator. A set in $Z$ is called {\em relatively compact} if its closure is compact. 
The operator $T$ is said to be {\em compact} if it is continuous and maps bounded sets into relatively compact sets.

\begin{lemma}[Schauder's fixed-point theorem~{\cite[Chapter 3]{zeidler1993nonlinear}}]
\label{fixedpoint}
If $\Z$ is a non-empty, closed, bounded, convex subset of $Z$ and $T: \Z \rightarrow \Z$ is a compact operator, then $T$ has a fixed-point.  
\end{lemma}

It is well-known that the relative compactness may be characterized by the uniform boundedness and equicontinuity. The next lemma states this.

\begin{lemma}[Arzel\`a–Ascoli theorem~{\cite[Theorem A5]{rudin91}}]
\label{ascoli}
A set $\Z \subset Z$ is relatively compact if 
and only if $\Z$ satisfies the following two conditions:
(a) uniform boundedness, i.e., $\sup_{f\in \Z} \|f\|_{\infty} < \infty$, and 
(b) equicontinuity, i.e., for every $\epsilon > 0$, there exists a constant $\delta >0 $ such that $\sup_{f \in \Z} \|f(\theta) - f(\theta')\| < \epsilon, \forall \theta, \theta'$ with $\|\theta - \theta'\| < \delta$.

\end{lemma}

\subsection{Relations between probability metrics}
In Section~\ref{se:stability}, we use several standard relations between probability metrics. 
For the reader's convenience, we recall below the comparison between the Kantorovich metric and the total variation metric established by Gibbs and Su~\cite[Theorem~4]{gibbs2002choosing}. 
This result is used to convert bounds in total variation distance into bounds in Kantorovich distance on bounded type spaces.
\begin{lemma}[{Gibbs and Su~\cite[Theorem~4]{gibbs2002choosing}}]
\label{metric}
Let $(\Xi,d)$ be a metric space with finite diameter
$\inmat{Diam}(\Xi):=\sup\{d(x,y):x,y\in \Xi\}$.
Then, for any $\mu,\nu\in\mathscr P(\Xi)$,
\begin{equation}
\label{eq:K-TV-upper}
\dd_{\mathrm K}(\mu,\nu)
\leq
\inmat{Diam}(\Xi)\,
\dd_{\mathrm{TV}}(\mu,\nu).
\end{equation}
If, in addition, $\Xi$ is finite and
$d_{\min}:=
\min_{x\neq y} d(x,y)>0,$
then
\begin{equation}
\label{eq:K-TV-lower}
d_{\min}\,
\dd_{\mathrm{TV}}(\mu,\nu)
\leq
\dd_{\mathrm K}(\mu,\nu).
\end{equation}
\end{lemma}
In our analysis, only \eqref{eq:K-TV-upper} is used. 
It allows us to control the Kantorovich discrepancy between conditional beliefs by their total variation discrepancy whenever the underlying type space is bounded. 
The lower bound \eqref{eq:K-TV-lower} is included for completeness.

\end{document}